\documentclass[sn-mathphys-num]{sn-jnl}

\usepackage{graphicx}%
\usepackage{multirow}%
\usepackage{amsmath,amssymb,amsfonts}%
\usepackage{amsthm}%
\usepackage{mathrsfs}%
\usepackage[title]{appendix}%
\usepackage{xcolor}%
\usepackage{textcomp}%
\usepackage{manyfoot}%
\usepackage{booktabs}%
\usepackage{algorithm}%
\usepackage{algorithmicx}%
\usepackage{algpseudocode}%
\usepackage{listings}
\usepackage{mathtools}   %extention of amsmath
\usepackage{physics}
\usepackage{stmaryrd}    %for symbols ex: \rrbracket
\usepackage{tikz-cd}     %for graphics/figures
\usepackage{esint}       %symbols integrals
\usepackage{subcaption}  %subfigures
\usepackage{orcidlink}
\usepackage[most]{tcolorbox}
\usepackage{aliascnt}
\usepackage{dsfont} %Id matrix

\usepackage{geometry}
\definecolor{burgundy}{rgb}{0.5, 0.0, 0.13}
\definecolor{cinereous}{rgb}{0.6, 0.51, 0.48}
\definecolor{lightsalmonpink}{rgb}{1.0, 0.6, 0.6}

\newtheoremstyle{thmstyleone}% Numbered
{18pt plus2pt minus1pt}% Space above
{18pt plus2pt minus1pt}% Space below
{\itshape}% Body font
{0pt}% Indent amount
{\bfseries}% Theorem head font
{}% Punctuation after theorem head
{.5em}% Space after theorem headi
{}% Theorem head spec (can be left empty, meaning `normal')
\theoremstyle{thmstyleone}
\numberwithin{equation}{section}
\newtheorem{theorem}{Theorem}[section]

\newaliascnt{corollary}{theorem}
\newtheorem{corollary}[corollary]{Corollary}
\aliascntresetthe{corollary}

\newaliascnt{lemma}{theorem}
\newtheorem{lemma}[lemma]{Lemma}
\aliascntresetthe{lemma}

\newaliascnt{definition}{theorem}
\newtheorem{definition}[definition]{Definition}
\aliascntresetthe{definition}

\newaliascnt{proposition}{theorem}
\newtheorem{proposition}[proposition]{Proposition}
\aliascntresetthe{proposition}

\newaliascnt{remark}{theorem}
\newtheorem{remark}[remark]{Remark}
\aliascntresetthe{remark}

\newaliascnt{notation}{theorem}
\newtheorem{notation}[notation]{Notation}
\aliascntresetthe{notation}

\newaliascnt{assumptions}{theorem}
\newtheorem{assumptions}[assumptions]{Assumptions}
\aliascntresetthe{assumptions}

\newaliascnt{example}{theorem}
\newtheorem{example}[example]{Example}
\aliascntresetthe{example}

\newcommand{\Eq}[2]{\begin{equation}\label{#1}\begin{aligned}#2 \end{aligned}\end{equation}}

\newcommand{\theo}[2]{\rbox{\begin{theorem}\label{#1} #2 \end{theorem}}}
\newcommand{\coro}[2]{\rbox{\begin{corollary}\label{#1} #2 \end{corollary}}}
\newcommand{\lem}[2]{\bbox{\begin{lemma}\label{#1} #2 \end{lemma}}}

\newcommand{\prop}[2]{\bbox{\begin{proposition}\label{#1} #2 \end{proposition}}}
\newcommand{\rem}[2]{\begin{remark}\label{#1} #2 \end{remark}}

\newcommand{\assum}[2]{\bbox{\begin{assumptions}\label{#1} #2 \end{assumptions}}}
\newcommand{\ex}[2]{\begin{example}\label{#1} #2 \end{example}}

\newcommand{\bbox}[1]{\begin{tcolorbox}[arc=0mm,oversize,colback=cinereous!3!white,colframe=cinereous!100!white]#1\end{tcolorbox}}
\newcommand{\rbox}[1]{\begin{tcolorbox}[arc=0mm,oversize,colback=purple!3!white,colframe=purple!100!white]#1\end{tcolorbox}}

\renewcommand{\geq}{\geqslant}
\renewcommand{\leq}{\leqslant}

\DeclareMathOperator{\sign}{sign}

\newcommand{\eins}{\mathds{1}}

\title{\textbf{Hard edge asymptotics of correlation functions between singular values and one eigenvalue}}

\author*[1,2]{\textsc{\fnm{Matthias} \sur{Allard}} \orcidlink{0000-0002-5682-424X}}\email{m.allard@unimelb.edu.au}

\affil[1]{\textit{\orgdiv{School of Mathematics and Statistics}, \orgname{University of Melbourne}, \orgaddress{\street{813 Swanston Street}, \city{Parkville, Melbourne}, \postcode{3010}, \state{Victoria}, \country{Australia}}}}
\affil[2]{\textit{\orgdiv{Department of Mathematics}, \orgname{KU Leuven}, \orgaddress{\street{Celestijnenlaan 200 B bus 2400}, \city{Leuven}, \postcode{3001},  \country{Belgium}}}}

\begin{document}

\abstract{Any square complex matrix of size $n\times n$ can be partially characterized by its $n$ eigenvalues and/or $n$ singular values. While no one-to-one correspondence exists between those two kinds of values on a deterministic level, for random complex matrices drawn from a bi-unitarily invariant ensemble, a bijection exists between the underlying singular value ensemble and the corresponding eigenvalue ensemble. This enabled the recent finding of an explicit formula for the joint probability density between $1$ eigenvalue and $k$ singular values, coined $1,k$-point function. We derive here the large $n$ asymptotic of the $1,k$-point function around the origin (hard edge) for a large subclass of bi-unitarily invariant ensembles called polynomial ensembles and its subclass Pólya ensembles. This latter subclass contains all Meijer-G ensembles and, in particular, Muttalib-Borodin ensembles and the classical Wishart-Laguerre (complex Ginibre), Jacobi (truncated unitary), Cauchy-Lorentz ensembles. We show that the latter three ensembles share the same asymptotic of the $1,k$-point function around the origin.
In the case of Jacobi ensembles, there exists another hard edge for the singular values, namely the upper edge of their support, which corresponds to a soft edge for the eigenvalue (soft-hard edge). We give the explicit large $n$ asymptotic of the $1,k$-point function around this soft-hard edge.}


% 60B20 RM proba
% 15B52 RM algebra

\maketitle
\tableofcontents

%-----------------------------------------------------------
\section{Introduction}\label{Introduction}
When studying a square complex matrix $X \in \mathbb{C}^{n\times n}$, one is usually interested in its \textit{eigenvalues} $z=\mathrm{diag}(z_1,\ldots,z_n)\in \mathbb{C}^n$, which are in general complex for non-Hermitian matrices, as they encode important information about the associated linear transformation. However, there exists another interesting set of values associated with the linear transformation: the \textit{singular values}. The squared singular values $a=\mathrm{diag}(a_1,\ldots,a_n)\in \mathbb{R}_+^n$ of $X$ are the eigenvalues of the matrix $X^\dagger X$, where $\dagger$ denotes the Hermitian conjugate. The singular values are then the positive square roots of the $a$'s.

While those two sets of values are often studied or used separately depending on the context, they can be complementary. Indeed, the information they capture is different because, albeit being related, there is no one-to-one correspondence between the two sets of values. This is why looking at both eigenvalues and singular values can be relevant and even crucial in certain areas. For instance, in the study of random non-Hermitian Hamiltonians and, in particular, in quantum chaos and open quantum systems \cite{Roccati2024,Porras2019,Herviou2019,Brunelli2023,Nandy2025}, the physics and the chaotic behavior of the system seems to be best captured by the singular values, without making the eigenvalues redundant \cite{Baggioli2025,Prasad2025}. The study of non-Hermitian quantum Hamiltonians via the singular values is sometimes referred to as biorthogonal quantum mechanics \cite{Brody2013}. The same complementarity is true in Quantum Chromodynamics~\cite{Kanazawa2011,Kanazawa2012}
 as well as in topological statistics of Hamiltonians~\cite{Braun2022,Hahn2023,Hahn2023a,Hajong2025} or even in Time Series Analysis of time-lagged matrices~\cite{Thurner2007,Long2023,Yao2022,Loubaton2021,Bhosale2018,Nowak2017}. On a practical level, computing a singular value decomposition is much more stable and faster than an eigenvalue decomposition, hence the use of both sets of values in numerical analysis \cite{An2024,Chen2024,Jones2024} to improve accuracy and algorithm complexity. 

 While, in general, for non-Hermitian matrices of fixed size $n$, there is no one-to-one correspondence between eigenvalues and singular values, there exists a unique equality between the two, coming from the determinant
 \Eq{eq:prodsev}{
|\det (X)|^2=|\det(z)|^2=\prod_{k=1}^n r_k=\det (X^\dagger X)=\det(a)=\prod_{k=1}^n a_k,
}
where $r_k$ are the squared eigenradii $r_k=|z_k|^2$. In addition, there exist various inequalities such as Weyl's inequalities~\cite{Weyl1949} from which it follows, in particular, that the largest singular value bounds the largest eigenradius from above and the smallest  eigenradius is bounded from below by the smallest singular value. These bounds might be the source of non-trivial correlations between eigenvalues and singular values which should survive in the limit of large matrix size.

On a probabilistic level, some bijection has been found \cite{Kieburg2016} between the probability densities of the singular values and eigenvalues for a large class of matrix ensembles: bi-unitarily invariant ensembles. If a density $f_{\rm BU}\in\mathrm{L}^1(\mathrm{GL}(n,\mathbb{C}))$ has the property
\Eq{}{
f_{\rm BU}(U_1 X U_2)=f_{\rm BU}(X)\quad{\rm for\ all}\  U_1,U_2 \in \mathrm{U}(n)\ {\rm and}\ X \in \mathrm{GL}(n,\mathbb{C}),
}
it is then called \textit{bi-unitarily invariant}. Many classical random matrix ensembles are bi-unitarily invariant. For instance, the induced Ginibre ensemble~\cite{Akemann2015} (also known as Laguerre, Wishart or chiral Gaussian unitary ensemble, when looking at the singular values)
\begin{equation}\label{eq:Lag buinv}
 f_{\rm Lag}(X)\propto \det (X^\dagger X)^\alpha\exp(-\Tr( X^\dagger X)),\ \alpha>-1,
\end{equation}
as well as the induced Jacobi ensemble~\cite{Akemann2015,Forrester2010} 
(also known as the ensemble of truncated unitary matrices)
\begin{equation}\label{eq:Jac buinv}
 f_{\rm Jac}(X)\propto\det (X^\dagger X)^\alpha \det (\eins_n-X^\dagger X)^\beta \Theta(\eins_n-X^\dagger X),\ \alpha>-1,\, \beta>0,
\end{equation}
with $\eins_n$ the $n\times n$ identity matrix and $\Theta$ the Heaviside step function equal to $1$ when the argument is a positive definite matrix and $0$ otherwise. Another example is the Cauchy-Lorentz ensemble~\cite{Wirtz2015}
\begin{equation}\label{eq:CL buinv}
 f_{\rm CL}(X)\propto\frac{\det (X^\dagger X)^\alpha }{\det (\eins_n+X^\dagger X)^{\beta+2(n-1)}},\ \alpha>-1,\, \beta-\alpha>1,
\end{equation}
and, more generally, any ensemble whose density depends only on $\det(P(X^\dagger X))$ and $\Tr( Q( X^\dagger X))$, where the functions $P,Q$ are any (convergent) series.

A direct consequence of the bijection between densities given in \cite{Kieburg2016} is that any bi-unitarily invariant ensemble is entirely determined by the density of its singular values and vice-versa; any density on $\mathrm{L}^1(\mathbb{R}_+^n)$ can be traced back to a bi-unitarily invariant density on $\mathrm{L}^1(\mathrm{GL}(n,\mathbb{C}))$. This is why we associate the ensemble on the singular values with the one on the full matrix space, e.g. the underlying ensemble on the singular values of the Ginibre ensemble \eqref{eq:Lag buinv} is the Laguerre ensemble.

In a recent article \cite{Allard2025b}, we have introduced the notion of $j,k$-point correlation functions between $j$ eigenvalues and $k$ singular values. For fixed matrix size $n$, we have exploited the bijection between the underlying densities of singular values and eigenvalues of a bi-unitarily invariant ensembles to get the $1,k$-point correlation functions between one eigenradius (modulus of the eigenvalue) and $k$ singular values. For $n>1$, the $1,k$-point correlation functions $f_{1,k}$ can be defined weakly, with the help of any continuous bounded test function $\phi \in C_b\left(\mathbb{R}_+^{1+k}\right)$, by the relation
\Eq{eq:1,k-point pdf}{
&\mathbb{E}\bigg[\frac{1}{n!} \sum_{\pi\in S_n} \phi(r(X); a_{\pi(1)}(X),\ldots,a_{\pi(k)}(X) ) \bigg]\\
=&\int_{\mathbb{R}_+^{1+k}}\phi(r;a_1,\ldots,a_k) f_{1,k}(r;a_1,\ldots,a_k)\,dr\, da_1\ldots da_k,  
}
where the expected value of a measurable function $\Phi:\mathrm{GL}(n,\mathbb{C})\to \mathbb{C}$ on $\mathrm{GL}(n,\mathbb{C})$ is explicitly defined by 
\Eq{eq:expected value}{
\mathbb{E}[\Phi]:=\int_{\mathrm{GL}(n,\mathbb{C})} \Phi(X)f_{\rm BU}(X) dX.
}
The $1,k$-point correlation functions are, thus, a generalization of the usual $k$-point correlation functions $R_k$, used to study determinantal point processes, up to a combinatorial factor~\cite{Akemann2015,Forrester2010} that makes them probability densities. The $1$-point probability function on the squared eigenradii and the $k$-point probability function on the squared singular values are respectively given by $f_{1,0}$ and $f_{0,k}$ which can be expressed as the marginals
\Eq{eq: 1 and k pt}{
f_{1,0}(r)&=\int_{\mathbb{R}_+^{k}}f_{1,k}(r;a_1,\ldots,a_k)\, da_1\ldots da_k,\\
f_{0,k}(a_1,\ldots,a_k)&=\int_{\mathbb{R}_+}f_{1,k}(r;a_1,\ldots,a_k)\,dr.
}
In \cite{Allard2025b}, we obtained a very general formula for the $1,k$-point correlation functions for any bi-unitarily invariant ensemble, which simplifies drastically when the induced probability density on the singular values is a \textit{polynomial ensemble}~\cite{Kuijlaars2014,Kuijlaars2014a,Kuijlaars2016,Foerster2020}, i.e. has the form 
\Eq{eq:polynomial ensemble def}{
f_{\rm SV}(x) =\frac{ \Delta_n(x)\det\left[ w_{k-1}(x_j)  \right]_{j,k=1}^n}{n!\,\det\left[ \mathcal{M}w_{k-1}(j)  \right]_{j,k=1}^n} ,
}
where $w_0,\ldots,w_{n-1}$ are weight functions on $\mathbb{R}_+$ such that $f_{\rm SV}$ is a probability density on $\mathbb{R}_+^n$. The $n$-dimensional Vandermonde determinant of an $n$-dimensional vector $x=(x_1,\ldots,x_n)$ is denoted by 
\begin{equation}\label{Vandermonde}
\Delta_n(x)=\det [x_j^{k-1}  ]_{j,k=1}^n=\prod_{1\leq j<k\leq n}(x_k-x_j).
\end{equation} We introduced, here, the Mellin transform  $\mathcal{M}$ on $\mathbb{R}_+$,
 \Eq{eq:M-trans}{
\mathcal{M}f(s):=\int_0^\infty dx\ x^{s-1}f(x)
}
for an $L^1(\mathbb{R}_+)$-function $f$ and $s\in\mathbb{C}$ such that the integral converges absolutely. Note that the structure of polynomial ensembles \eqref{eq:polynomial ensemble def} appears very naturally in random matrices as the Jacobian of the change of coordinates to go on the singular values yields a squared Vandermonde determinant.

The results in \cite{Allard2025b} are made even more explicit for a certain type of polynomial ensemble,  namely the \textit{P\'olya ensembles} (formerly coined polynomial ensembles of derivative type)~\cite{Kieburg2015,Kieburg2016,Kieburg2022,Foerster2020}. A polynomial ensemble is a Pólya ensemble if there exists $w$ such that 
 \Eq{eq: polya w}{
 w_{j}(x)= (-x\partial_x)^{j} w(x) \in \mathrm{L}^1(\mathbb{R}_+),\qquad \forall j \in \llbracket 0,n-1\rrbracket.
}
To guarantee that we deal with probability measures it has been shown in~\cite{Foerster2020} that $w$ is then related to Pólya frequency functions. Many classical ensembles are Pólya ensembles. It is, for instance, the case of the Laguerre \eqref{eq:Lag buinv}, Jacobi \eqref{eq:Jac buinv} and Cauchy-Lorentz \eqref{eq:CL buinv} ensembles, whose respective corresponding Pólya weight functions \eqref{eq: polya w} are
\Eq{eq: polya weights}{
w_{\rm Lag}(x)=x^\alpha e^{-x},\quad w_{\rm Jac}(x)=x^\alpha (1-x)^{\beta+n-1} \Theta(1-x),\quad w_{\rm CL}(x)=\frac{x^\alpha}{(1+x)^{\beta+n-1}}.
}
Those ensembles are well studied and appear in many different applications such as the study of wireless communication systems, entanglement of a random pure quantum state, quantum conductance; see \cite{Forrester2010,Akemann2015,LalMehta2004}. Note that the Cauchy-Lorentz ensemble is related to the circular unitary ensemble (CUE)---and its truncated version---via a Cayley transformation. On the level of eigenvalues it corresponds to a simple change of variable, namely a stereographic projection from the unit circle to the real line (see \cite[p.68]{Forrester2010}).

Other important ensembles falling in the class of Pólya ensemble are the Muttalib-Borodin ensembles \cite{Borodin1998,Forrester2017} which have been introduced to model two-body interactions in random matrix models \cite{Muttalib1995} and used, for instance, in the study of disordered conductors, quantum transport, Brownian motion \cite{Beenakker1997,Lueck2006,Takahashi2012}. We give the example of the Laguerre type Muttalib-Borodin ensemble whose Pólya weight is given by
\Eq{eq: polya MB}{
w_{\rm MB}(x)=x^\alpha e^{-x^\theta},\quad \alpha>-1, \theta>0.
}

In fact, not only are most of the classical random matrix ensembles Pólya ensembles but also their composition, in the sense that given two complex square matrices $X$ and $Y$ drawn from different (or same) Pólya ensembles, then the matrix $XY$ also follows a Pólya ensemble, whose corresponding Pólya weight is given by the multiplicative convolution of the two original weights, cf. \cite{Foerster2020,Kieburg2019,Kuijlaars2014,Kuijlaars2016}. This nice closure property makes the study of this class of ensemble all the more interesting.

Coming back to general polynomial ensembles, they enjoy some additional structure from their belonging to a much larger class of ensembles called \textit{determinantal point processes}; see~\cite{Anderson2010,Akemann2015,Deift2009,Forrester2010}. This means that the joint probability distribution $f_{\rm SV}$ can be written in the form
\Eq{eq: DPP poly ens}{
f_{\rm SV}(x) =\frac{1}{n!}\det\left[K_n(x_j,x_k)  \right]_{j,k=1}^n,
}
where $K_n$ is the kernel function, and all the $k$-point correlation functions have a similar form where only the size of the determinant changes. Note that $K_n$ is not uniquely given. Indeed, due to elementary properties of the determinant, for a non-vanishing function $g$, the kernel $[g(x_1)/g(x_2)]K_n(x_1,x_2)$ is also a correlation kernel for the same point process. However, for polynomial ensembles one can choose $K_n$ to be polynomial of degree $n-1$ in the second entry, which thus makes it unique. It is important to stress that not only does the kernel play a crucial role in the point process of the singular values but also in the point process of the eigenradii and, a fortiori, in the combined point process as it can be seen in the main results of \cite{Allard2025b}, reminded here for convenience.

 %=============THEO POLY ENSEMBLE==============
\theo{theo:poly ensemble}{\textup{(Theorem 1.4 \cite{Allard2025b})} Let $n\in \mathbb{N}$, $n>2$, $k\in \llbracket 1,n\rrbracket$ and consider a random matrix drawn from a bi-unitarily invariant ensemble on $\mathrm{GL}(n,\mathbb{C})$  having a polynomial ensemble with joint probability density~\eqref{eq: DPP poly ens} for the squared singular values. The $1,k$-point correlation function between one squared eigenradius and $k$ squared singular values is given by
\Eq{eq:1,kpt poly}{
f_{1,k}(r; a_1,\ldots,a_k) =&\frac{(n-k)!}{(n-1)!}\int_{0}^{\infty} dt \int_0^r\frac{dv}{v} \varphi_n\left(\tfrac{v}{r},t\right)\\
&\times\det\left(\begin{array}{c c} 
    	 K_n\left(v,-rt\right)\quad & K_n\left(v,a_c\right)-\delta(v-a_c)  \\
      K_n\left(a_b,-rt\right) \quad & K_n(a_b,a_c) 
\end{array}\right)_{b,c=1}^k.
}
with $\delta$ the Dirac delta function,
\Eq{eq: def phi}{
 \varphi_n(x,t):= x(1-x)^{n-2}(1+t)^{-(n+2)}\left[\left(1-\frac{x}{n}\right)(1+t)-(1-x)\left(1+\frac{1}{n}\right) \right]
}
with $K_n$, the correlation kernel of the polynomial ensemble, chosen to be a polynomial of degree $n-1$ in its second argument. 

The $1$-point functions, respectively on one squared singular value and one squared eigenradius, are given by
\Eq{}{\rho_{\rm SV}(a):=f_{0,1}(a)=\frac{1}{n}K_n(a,a)
}
and
\Eq{eq:1pointreal2}{
\rho_{\rm EV}(r):=f_{1,0}(r)=n \int_{0}^{\infty} dt \int_0^r \frac{dv}{v} \varphi_n\left(\frac{v}{r},t\right) K_n\left(v,-rt\right).}
}
While the $1,k$-point correlation functions contains all the information about the correlation between one eigenradius and $k$ singular values, a way to extract this information consists in looking at the difference between the $1,k$-point function $f_{1,k}$ and the product of the respective $1$-point function  $f_{1,0}$ on the eigenradii and the $k$-point function $f_{0,k}$ on the singular values \eqref{eq: 1 and k pt}. This is then a measurement of the distance to statistical independence: a covariance measurement. The resulting function has been coined \textit{$1,k$-cross-covariance density} in \cite{Allard2025b} and is denoted by
\Eq{eq: def 1,k cov}{
\mathrm{cov}_{1,k}(r;a_1,\ldots,a_k):=f_{1,k}(r;a_1,\ldots,a_k)-f_{1,0}(r)f_{0,k}(a_1,\ldots,a_k).
}
As a shorthand, we will refer to the $1,1$-cross-covariance density, simply, as \textit{cross-covariance density} and denote it $\mathrm{cov}:=\mathrm{cov}_{1,1}$. It was then shown in \cite{Allard2025b} that, for Pólya ensembles, the $1,k$-cross-covariance density could be put in the form given by the following proposition.
%=============CORO POLY ENSEMBLE EXPLICIT==============

\prop{prop:poly ensemble explicit}{\textup{(Prop.1.7 \cite{Allard2025b})} Let $n\in \mathbb{N}$, $n>2$.  With the same assumptions and notations as in \autoref{theo:poly ensemble}, we assume that $f_{\rm SV}$ is the joint probability density of the squared singular values of a Pólya ensemble associated to an $n$-times differentiable weight function $w_n\in C^n(\mathbb{R}_+)$. Then, the $1,k$-cross-covariance density function has the form
\Eq{cov.prop}{
\mathrm{cov}_{1,k}(r;a_1,\ldots,a_k)= \frac{(n-k)!}{(n-1)!}\left.\partial_\mu\det[K_n(a_b,a_c)+\mu\, \mathbf{C}_n(r;a_b,a_c)]_{b,c=1}^k\right|_{\mu=0}
}
with
\Eq{cov.prop2}{
\mathbf{C}_n(r;a_1,a_2)=& \sum_{\gamma=0,1} H_\gamma(r,a_2)\left[\Theta(r-a_1) \frac{1}{a_1} \Psi_\gamma\left(\frac{a_1}{r}\right)- V_\gamma(r,a_1) \right]
}
with $\Theta$ the Heaviside step function and for $\gamma=0,1$ we have employed the functions
\begin{eqnarray}
\Psi_\gamma(x)&:=& \left(\frac{1-x}{nx-1}\right)^\gamma x(1-x)^{n-2} \left(nx-1\right),\label{psi.def}\\
H_\gamma(x,y)& :=&\int_0^1 du\ q_n(yu) \partial_{u}^\gamma \left[ u^\gamma\label{H.def} \frac{\rho_{\rm EV}(xu)}{w_n(xu)}\right] , \\
V_\gamma(x,y) &=&  \int_0^1 du\  p_{n-1}\left(yu\right) \left(u \partial_{u}\right)^{1-\gamma} w_n(xu),\label{V.def}
\end{eqnarray}
where $p_{n-1}$ and $q_n$ are the bi-orthonormal pair of functions composing the kernel~\eqref{eq: Kernel poly ens} of  $f_{\rm SV}$ which can be expressed as
\Eq{eq:biorthogonal}{
p_{n-1}(x) &=\sum_{c=0}^{n-1} \binom{n-1}{c} \frac{(-x)^c}{\mathcal{M}w_n(c+1)} \quad{\rm and}\quad q_n(x) = \frac{1}{n!}\partial_x^n[x^n w_n(x)],
}
according to~\cite[Lemma 4.2]{Kieburg2016}.
}

The goal of this article is to exploit the results of \cite{Allard2025b} where everything has been done for fixed matrix size $n$ and study, now, the large $n$ limit.  Let us stress that all the results of this article are completely new as the results of \cite{Allard2025b} were new.

The article is organized as follows. First, in \autoref{Double scaling limit at the origin}, we look at the double scaling limit of both the $1,k$-point correlation function and the $1,k$-cross-covariance density around the origin for polynomial ensembles. This is interesting and non trivial, as most polynomial ensembles will have a hard edge at the origin. Indeed, unless there is a strong repulsion from the origin, the lower edge of the support will be a hard wall due to the definition of eigenradii and singular values, which are non negative. Moreover, the deterministic  Weyl's inequalities~\cite{Weyl1949} between eigenradii and singular values imply non trivial correlations between the two, where they share the same hard edge. The main results of this section are \autoref{theo: 1,k poly ensemble}, \autoref{theo: 1,k polya ensemble} and \autoref{coro:polya ensemble}. We apply the results to the example of the Muttalib-Borodin ensembles.

Some ensembles can have another hard edge, as it is the case of Jacobi ensembles. In \autoref{Soft-hard edge scaling limit for Jacobi ensembles}, we look, in particular, at the Jacobi ensemble  \eqref{eq:Jac buinv} around the upper hard edge for the singular values, which corresponds to a soft edge for the eigenradius. The main result of this section is given by \autoref{theo:jac up HE}.
In \autoref{Proof HE} and \autoref{Proof SH}, we give the proofs for, respectively, the hard edge at the origin and the soft-hard edge. Then, we discuss the results and give some plots in \autoref{Discussion}. Additional computation regarding the consistency of our assumptions are given in \autoref{Appendix}.

% %-----------------------------------------------------------
\section{Double scaling limit at the origin}\label{Double scaling limit at the origin}
\subsection{Polynomial ensembles}
Let us start with some heuristics. When looking at \autoref{theo:poly ensemble}, the double scaling limit (if existent) of the $1,k$-point correlation function will be given by scaling appropriately the squared eigenradius and the squared singular values and then taking the large $n$ limit. One would expect the result to involve the double scaling limit of the kernels $K_n$ and the double scaling limit of the $\varphi_n$. 

The goal is then to find the correct scalings. To begin, we examine \eqref{eq:1,kpt poly} and first proceed with the change of variable $v\mapsto rv$. The arguments of $\varphi_n$ do not depend on $r$ or $a$ anymore and one can see that the term $n^2 \varphi_n\left(x/n,t/n \right)$ admits a pointwise limit
\Eq{}{
\lim_{n\to\infty}n^2 \varphi_n\left(\frac{x}{n},\frac{t}{n} \right)=x(t+x-1)e^{-x} e^{-t}.
}
This hints one needs to proceed with the change of integration variables $(v,t)\mapsto (v/n,t/n)$.

Now, let us assume that the appropriate scaling for the squared eigenradius is given by $\nu_n$, thus we let $r\mapsto r/\nu_n$. It becomes clear that if the result involves the double scaling limit of the kernels, the correct scaling should be
\Eq{}{
\frac{1}{n\nu_n}K_n\left(\frac{x}{n\nu_n},\frac{y}{n\nu_n}\right).
}
To be consistent, the correct scaling for the squared singular values is then $n\nu_n$, thus we set $(a_1,\ldots,a_k)\mapsto \left(\frac{a_1}{n\nu_n},\ldots, \frac{a_k}{n\nu_n}\right)$. This surprising difference between the scaling of the squared eigenradius and squared singular values is discussed later.

Accounting for the change of measure and multiplying by the appropriate global scaling, the correct scaling of the $1,k$-point correlation function between one squared eigenradius and $k$ squared singular values should be
\Eq{eq: scaled 1,k pt2}{
&\frac{n^{k+1}}{\nu_n (n\nu_n)^{k}}f_{1,k}\left(\frac{r}{\nu_n}; \frac{a_1}{n\nu_n},\ldots, \frac{a_k}{n\nu_n}\right)\\
&\quad =\frac{n^k (n-k)!}{n!}\int_{0}^{\infty} dt \int_0^{n}\frac{dv}{v} n^2\varphi_n\left(\frac{v}{n }, \frac{t}{n}\right)\\
&\quad \ \times \det\left(\begin{array}{c c} 
    	 \frac{1}{n\nu_n}K_n\left(\frac{rv}{n\nu_n},\frac{-rt}{n\nu_n}\right)\quad & \frac{1}{n\nu_n}K_n\left(\frac{rv}{n\nu_n},\frac{a_c}{n\nu_n}\right)-\delta(rv-a_c)  \\
     \frac{1}{n\nu_n} K_n\left(\frac{a_b}{n\nu_n},\frac{-rt}{n\nu_n}\right) \quad & \frac{1}{n\nu_n}K_n\left(\frac{a_b}{n\nu_n},\frac{a_c}{n\nu_n}\right) 
\end{array}\right)_{b,c=1}^k,
}
where we have used the scaling property of the Dirac delta function
\Eq{}{
\forall c\neq 0, \quad \delta\left(\frac{x}{c}\right)=c \,\delta(x).
}
To be allowed to take the limit $n\to \infty$ on both sides of \eqref{eq: scaled 1,k pt2} and then pass the limit below the integrals in the RHS, all limits must exist and one needs to justify exchanging limit and integrals. This is why one needs to make some restrictions on the kernel of the polynomial ensemble as its double scaling limit is not necessarily a function anymore or simply does not exist. 
We will make the following assumptions to guarantee existence of the limits and integrability of the different integrands involved in the computations. 

Let us first recall that the kernel $K_n$ of any polynomial ensemble \eqref{eq: DPP poly ens} can be cast in the form
\Eq{eq: Kernel poly ens}{
K_n(x,y)=\sum_{j=0}^{n-1}q_j(x)p_j(y),
}
with $p_j$ a polynomial of degree $j$ and $q_{j}$ chosen such that the two sets of functions $\{q_j\}_{j=0,\ldots,n-1}$ and  $\{p_j\}_{j=0,\ldots,n-1}$ satisfy the bi-orthonormality relation
\Eq{eq: bi-orthonormality}{
\int_{0}^\infty q_j(x)p_k(x)dx=\delta_{j,k},
}
where $\delta_{j,k}$ is the Kronecker delta function. Note that this does not uniquely define the functions $p_j$ and $q_j$ as one can always decompose the polynomials $p_j$ in another basis of polynomials up to degree $n-1$ and then take the corresponding linear combination of $q_j$ to form a new bi-orthonormal system. In particular, given a sequence $\{c_j\}_{j=0,\ldots,n-1}$ of non zero real numbers, $\{c_j q_j\}_{j=0,\ldots,n-1}$ and  $\{c_j^{-1}p_j\}_{j=0,\ldots,n-1}$ are also possible bi-orthonormal families. Thus, the $p_j$ are not canonically associated to the ensemble. Note also that $p_j$, $q_j$ might depend on $n$, we, however, do not make the notation transparent. 

%---assum:exist HE PE--------------------------------------------------------------------------
\assum{assum:exist HE PE}{Let $n\in \mathbb{N}$. Let $\{p_j\}_{j=0,\ldots,n-1}$, $\{q_j\}_{j=0,\ldots,n-1}$ be bi-orthonormal families of functions composing the kernel $K_n$ of a polynomial ensemble as defined in \eqref{eq: Kernel poly ens}. 
\begin{enumerate}
\item There exists a positive sequence $\{\nu_n\}_{n\in\mathbb{N}}$, such that for $(x,y)\in\mathbb{R}_+\times \mathbb{R}$, the point-wise limit
\Eq{eq: def limit kernel}{
K^{(\infty)}(x,y):=\lim_{n\to \infty} \frac{1}{n \nu_n}K_n\left(\frac{x}{n \nu_n},\frac{y}{n \nu_n}\right)
}
is a measurable function on $\mathbb{R}_+\times \mathbb{R}$.

\item There exist two real continuous functions $f$ and $g$, independent of $n$ such that:
\begin{enumerate}
\item For almost all $x$ in $\mathbb{R}_+$,
\Eq{eq: bound q_j}{
\forall j\in\llbracket 0,n-1\rrbracket,  \quad \left|\frac{1}{\nu_n}q_j\left(\frac{x}{n \nu_n}\right)\right|\leq \exp(f(x)),
}
with $f(u)=o(u)$, as $u\to \infty$, and $u\mapsto\exp(f(u))\in \mathrm{L}^1_{\rm loc}(\mathbb{R}_+)$ locally integrable.
\item For almost all $y$ in $\mathbb{R}$,
\Eq{eq: bound p_j}{
\forall j\in\llbracket 0,n-1\rrbracket,  \quad \left|p_j\left(\frac{y}{n \nu_n}\right)\right|\leq \exp(g(y)).
}
\end{enumerate}
\item We have, for $R>0$,
\Eq{eq: asymp pj polynom}{
\lim_{R\to \infty} \lim_{n\to \infty}\int_{R}^{\infty} dt  \frac{ (1+t)}{\left(1+\frac{t}{n}\right)^{n+2}}\frac{1}{n}\sum_{j=0}^{n-1} \abs{p_j\left(\frac{x t}{n\nu_n} \right) }=0.
}
\end{enumerate}
}
\rem{}{
The bounds \eqref{eq: bound p_j} and \eqref{eq: bound q_j} are not too restrictive in regard to the assumption on the existence of the limiting kernel \eqref{eq: def limit kernel}. Indeed, if a limiting kernel exists around the origin, then it is reasonable to assume that the products $q_j(x)p_j(y)$, remains bounded independently of $n$ on the scale $1/(n\nu_n)$ around the origin. Similarly, while \eqref{eq: asymp pj polynom} is only technical, it guarantees that the contribution of the polynomials to the integral in \eqref{eq: scaled 1,k pt2} remains in a bounded interval near $0$, i.e. near the hard edge.
}
\rem{}{
While finding a set of bi-orthonormal functions $p_j$, $q_j$ is not easy in general, once found for any $n$, the conditions are rather easy to check. Note that imposing the conditions \eqref{eq: bound p_j} and \eqref{eq: bound q_j} separately on $p_j$ and $q_j$ constrains the choice of the bi-orthonormal system without making it unique. 
}
\ex{}{
For Jacobi, Cauchy-Lorentz and Laguerre ensembles, $K^{(\infty)}$ is expressed in terms of the Bessel kernel, cf. \autoref{prop: limit kernel classic Pólya}, and a possible choice for the functions $p_j$ and $q_j$ is given by \eqref{eq:polya ortho functions} which one needs to rescale by a factor $\xi_n$ and consider instead the families $\{\xi_n^{-1} p_j\}_{j=0,\ldots,n-1}$, $\{\xi_n q_j\}_{j=0,\ldots,n-1}$. We take $\nu_n=1$ and $\xi_n=1$, for Laguerre, and $\nu_n=n$ and $\xi_n=n^{1+\alpha}$, for Jacobi and Cauchy-Lorentz.

For polynomials $p_j$ one can easily find the bound, for $j=0,\ldots,n-1$,
\Eq{}{
\left|\frac{1}{\xi_n}p_j\left(\frac{y}{n \nu_n}\right)\right|\leq \exp(\abs{y}+c), \quad c>0,
}
and take $g(y)=\abs{y}+c$; see \eqref{eq: bound pj polya}.

Then, as long as the P\'olya weights remains bounded on $\mathbb{R}_+$, which corresponds to $\alpha\geq 0$ for these ensembles, the corresponding functions $q_j$ are simply bounded by a $n$-independent constant
\Eq{}{
\left|\frac{\xi_n}{\nu_n}q_j\left(\frac{x}{n \nu_n}\right)\right|\leq C,\quad C>0,
}
for $j=0,\ldots,n-1$. Thus one can take $f(x)=\ln(C)$.

In the case the P\'olya weight is unbounded, i.e. $\alpha<0$, one must either select the biorthonormal functions more carefully---together with the general bound \cite[Corollary 2.6]{Kieburg2022}---or exploit the integrable structure of the ensembles, specifically the Christoffel-Darboux formula or the integral representation of the kernel \eqref{eq: integ kernel polya}, to verify that it is indeed uniformly bounded with respect to $n$.

The last assumption \eqref{eq: asymp pj polynom} requires more work to check and it is done in \autoref{Appendix}.
}
Under these assumptions, the double scaling limit version of \autoref{theo:poly ensemble} around the origin, in the large $n$ limit, yields \autoref{theo: 1,k poly ensemble}. Let us stress, here, that the eigenradius and the singular values scale differently. The scale of the smallest squared eigenradius is $n$ times larger than the scale of the smallest squared singular value. This was already observed for classical ensembles like the Laguerre ensemble \eqref{eq:Lag buinv}. Indeed, in this case the limiting macroscopic distribution of the singular values is the quarter circle law and the mean level spacing around the origin is of order $1/n$. As for the eigenvalues, they tend to the uniform distribution in a annulus/disk. Thus, the smallest eigenvalue scales as $1/\sqrt{n}$ and so is the smallest eigenradius. Going in squared coordinates yields this scaling ratio of $1/n$. A general and simple explanation of why this scaling ratio is universal is lacking. However, from the analysis perspective, this ratio appears very naturally and is indeed universal, at least for polynomial ensembles, as shown in \autoref{theo: 1,k poly ensemble}.

%=============THEO POLY ENSEMBLE==============
\theo{theo: 1,k poly ensemble}{
Let $k\in \mathbb{N}$ and consider a random matrix drawn from a bi-unitarily invariant ensemble on $\mathrm{GL}(n,\mathbb{C})$  having a polynomial ensemble with joint probability density~\eqref{eq: DPP poly ens} for the squared singular values. Under \autoref{assum:exist HE PE}, the point-wise limit
\Eq{}{f_{1,k}^{(\infty)}(r; a_1,\ldots,a_k):=\lim_{n\to \infty}\frac{n}{\nu_n^{k+1}}f_{1,k}\left(\frac{r}{\nu_n}; \frac{a_1}{n\nu_n},\ldots, \frac{a_k}{n\nu_n}\right)}
of the $1,k$-point correlation function between one squared eigenradius and $k$ squared singular values is given by
\Eq{eq:1,kpt poly inf}{
f_{1,k}^{(\infty)}(r; a_1,\ldots,a_k)=&
\int_{0}^{\infty} dt  \int_0^\infty \frac{dv}{v} \varphi^{(\infty)}\left(\tfrac{v}{r},t\right)\\
&\times\det\left(\begin{array}{c c} 
    	 K^{(\infty)}\left(v,-rt\right)\quad & K^{(\infty)}\left(v,a_c\right)-\delta(v-a_c)  \\
      K^{(\infty)}\left(a_b,-rt\right) \quad & K^{(\infty)}(a_b,a_c) 
\end{array}\right)_{b,c=1}^k,
}
with $\delta$ the Dirac delta function and
\Eq{eq: def phi inf}{
 \varphi^{(\infty)}(x,t):= x(t+x-1)e^{-x} e^{-t}.
}
The point-wise limit of $1$-point functions, respectively on one squared singular value and one squared eigenradius, are given by
\Eq{eq:1pointreal SV}{\rho_{SV}^{(\infty)}(a):=\lim_{n\to\infty} \tfrac{1}{\nu_n}\rho_{SV}\left(\tfrac{a}{n\nu_n} \right)=K^{(\infty)}(a,a)
}
and
\Eq{eq:1pointreal}{
\rho_{EV}^{(\infty)}(r):=\lim_{n\to\infty} \tfrac{n}{\nu_n}\rho_{EV}\left(\tfrac{r}{\nu_n} \right)= \int_{0}^{\infty}dt  \int_0^\infty \frac{dv}{v} \varphi^{(\infty)}\left(\tfrac{v}{r},t\right) K^{(\infty)}\left(v,-r t\right).
}
}
\rem{rem: speed decrease 1}{
Rewriting the result as follows
\Eq{eq: order scaling}{
\frac{1}{\nu_n (n\nu_n)^{k}}f_{1,k}\left(\frac{r}{\nu_n}; \frac{a_1}{n\nu_n},\ldots, \frac{a_k}{n\nu_n}\right)\underset{n\to \infty}{=}\frac{1}{n^{1+k}}\left[f_{1,k}^{(\infty)}(r; a_1,\ldots,a_k)+o(1)\right],
}
one can see that the scaling limit of the $1,k$-point function at the origin is of order $O(1/n^{1+k})$ as $n\to \infty$. We stress that \eqref{eq: order scaling} is on the level of probability densities, which are thus properly normalized. It should be noted, however, that $\rho_{EV}^{(\infty)},\rho_{SV}, ^{(\infty)}$ and $f_{1,k}^{(\infty)}$ are not probability densities.
}

The double scaling limit of the $1,k$-cross-covariance density follows directly from the definition \eqref{eq: def 1,k cov} and is given by the following corollary. 

%=============THEOREM POLY ENSEMBLE==============
\coro{coro: cov poly ensemble}{ 
With the same assumptions as in \autoref{theo: 1,k poly ensemble}, the point-wise limit
\Eq{}{
\mathrm{cov}_{1,k}^{(\infty)}(r; a_1,\ldots,a_k):=\lim_{n\to \infty}\frac{n}{\nu_n^{k+1}}\mathrm{cov}_{1,k}\left(\frac{r}{\nu_n}; \frac{a_1}{n\nu_n},\ldots, \frac{a_k}{n\nu_n}\right)
}
of the $1,k$-cross-covariance density functions between one squared eigenradius and $k$ squared singular values is given by
\Eq{cov.prop inf}{
\mathrm{cov}_{1,k}^{(\infty)}(r; a_1,\ldots,a_k) = \partial_\mu\det\left[K^{(\infty)}(a_b,a_c)+\mu\, \mathbf{C}^{(\infty)}(r;a_b,a_c)\right]_{b,c=1}^k\bigg|_{\mu=0}
}
with
\Eq{eq: Cinf integ}{
\mathbf{C}^{(\infty)}(r;a_b,a_c):=\int _{0}^{\infty}dt  \int_0^\infty \frac{dv}{v} \varphi^{(\infty)}\left(\tfrac{v}{r},t\right) K^{(\infty)}\left(a_b,-rt\right)\left[\delta(v-a_c)-K^{(\infty)}(v,a_c)\right],
}
where $\delta$ is the Dirac delta function. In particular for $k=1$, the limiting cross-covariance reads
\Eq{}{
  \mathrm{cov}^{(\infty)}(r; a)&=\mathbf{C}^{(\infty)}(r;a,a).
}
}
\begin{proof}
Starting from \autoref{theo: 1,k poly ensemble} and using the definition \eqref{eq: def 1,k cov}, the scaling limit of the $1,k$-cross-covariance density is given explicitly by
\Eq{}{
\mathrm{cov}_{1,k}^{(\infty)}(r;a_1,\ldots,a_k)=&\int_{0}^{\infty} dt  \int_0^\infty \frac{dv}{v} \varphi^{(\infty)}\left(\tfrac{v}{r},t\right)\\
&\times\det\left(\begin{array}{c c} 
    	 0\quad & K^{(\infty)}\left(v,a_c\right)-\delta(v-a_c)  \\
      K^{(\infty)}\left(a_b,-rt\right) \quad & K^{(\infty)}(a_b,a_c) 
\end{array}\right)_{b,c=1}^k\\
=&\sum_{l,m=1}^k(-1)^{l+m}\det[K^{(\infty)}(a_b,a_c)]_{\substack{b\in\llbracket 1,k\rrbracket\setminus\{l\}\\c\in\llbracket 1,k\rrbracket\setminus\{m\}}}\mathbf{C}^{(\infty)}(r;a_l,a_m).
}
Applying Jacobi's formula for derivative of determinants yields the result.
\end{proof}
%---------------------------------------------------------------------------
\subsection{Pólya ensembles}
For the subclass of Pólya ensembles, the formulas of \autoref{theo: 1,k poly ensemble} can be made more explicit by computing the double integral. Furthermore, as Pólya ensembles enjoy more structure and are entirely defined by their weight function, one can conveniently translate and make precise \autoref{assum:exist HE PE} directly on the Pólya weight. 

Thus, in the case of Pólya ensembles, the following assumptions (cf. \cite[Assumptions 2.2]{Kieburg2022}) are sufficient for the existence of a hard edge at the origin and guarantee the existence of the limiting kernel \eqref{eq: def limit kernel}. The final assumption is an additional condition not present in \cite[Assumptions 2.2]{Kieburg2022}, whereas the remaining assumptions have been adapted to accommodate the new scalings $\xi_n$ and $\nu_n$.

%---assum:exist HE--------------------------------------------------------------------------
\assum{assum:exist HE}{Let $n\in \mathbb{N}$. We assume that the analytic continuation of the Mellin transform $\mathcal{M}w_n$ is holomorphic on $\mathbb{C} \setminus \left(] - \infty, 1[\cup]n, \infty[\right)$ and that there exist two real positive sequences, $\{\xi_n\}_{n\in\mathbb{N}}$ and $\{\nu_n\}_{n\in\mathbb{N}}$, such that:
\begin{enumerate}
\item For all fixed $x\in \mathbb{R}_+$ and fixed $s\in \mathbb{C}\setminus]-\infty,1]$, the point-wise limits
\Eq{eq: def limit polya w}{
\lim_{n\to \infty} \frac{\xi_n}{\nu_n}w_n\left(\frac{x}{\nu_n} \right):=w^{(\infty)}(x),\quad \lim_{n\to \infty} \xi_n\nu_n^{s-1}\mathcal{M}w_n(s):=\Tilde{w}^{(\infty)}(s)
}
exist and are measurable functions. Then, $\Tilde{w}^{(\infty)}(s)=\mathcal{M}w^{(\infty)}(s)$.
\item There exists a constant $C>0$, $n$-independent, such that 
\Eq{eq: bound Tildew}{
\frac{1}{\xi_n\nu_n^{s-1}\abs{\mathcal{M}w_n(s)}}\leq C, \quad  s\in [1,n].
}
\item We have the domination
\Eq{eq: dom Tildew}{
\sup_{\arg(z)=\vartheta}\{ |\xi_n\nu_n^z \mathcal{M}w_n(1+z)|\}\leq D(\vartheta), \quad \forall n\in \mathbb{N},
}
where $0<D(\vartheta)<\infty$, for all $\vartheta \in [\pi/2,\pi[$.
\item There exists $j_*\in \llbracket 0,n-1\rrbracket$, and a positive function $h$, both $n$-independent, such that 
\Eq{eq: lower bound h}{
\frac{1}{\xi_n\nu_n^j \abs{\mathcal{M}w_n(j+1)}}\leq\frac{1}{h(j)}\leq C, \quad j\in[j_*,n-1],
}
and the function $h$ has the root asymptotic
\Eq{eq: root asymp bound h}{
\limsup_{j\to\infty} \frac{1}{h(j)^{1/j}}=0.
}
\end{enumerate} 
}

\rem{}{ The scaling $\xi_n$ does not explicitly play a role in the further computation. To see this one can look at \eqref{eq: integ scaled kernel polya}, \autoref{prop: kernel HE polya} and notice that only ratios $\mathcal{M}w_n(x)/\mathcal{M}w_n(y)$ are involved. This is reminiscent of the fact the choice of bi-orthonormal functions composing the kernel is not uniquely given.}
\ex{ex:classic polya}{For Laguerre, Jacobi and Cauchy-Lorentz ensembles, the respective Pólya weights are given by \eqref{eq: polya weights} where $\alpha$, $\beta$ are fixed, i.e. independent of $n$. The Laguerre weight being independent of $n$, there is no need to rescale, thus $\nu_n=1$ and $\xi_n=1$. For Jacobi and Cauchy-Lorentz weights, taking $\nu_n=n$ and $\xi_n=n^{1+\alpha}$ yields
\Eq{}{
w_{\rm Lag}^{(\infty)}(x)=w_{\rm Jac}^{(\infty)}(x)=w_{\rm CL}^{(\infty)}(x)=x^\alpha e^{-x}.
}
Starting from the Mellin transform of the respective weights
\Eq{}{
&\mathcal{M}w_{\rm Lag}(s)=\Gamma(s+\alpha),\quad \mathcal{M}w_{\rm Jac}(s)=\mathrm{B}(s+\alpha,\beta+n),\\
&\mathcal{M}w_{\rm CL}(s)=\mathrm{B}(s+\alpha,n-1-s+\beta-\alpha),
}
it is easy to check that we have, as expected,
\Eq{eq: Mellin polya weights}{
\Tilde{w}_{\rm Lag}^{(\infty)}(s)=\Tilde{w}_{\rm Jac}^{(\infty)}(s)=\Tilde{w}_{\rm CL}^{(\infty)}(s)=\Gamma(s+\alpha).
}
Regarding the $n$-independent bound $h$ \eqref{eq: lower bound h}, for Laguerre, Jacobi and Cauchy-Lorentz ensembles we have, respectively,
\Eq{}{
\xi_n\nu_n^{j}\mathcal{M}w_{\rm Lag}(j+1)&=\Gamma(j+1+\alpha),\quad \alpha>-1,\\
\xi_n\nu_n^{j}\mathcal{M}w_{\rm Jac}(j+1)&=n^{j+1+\alpha}\mathrm{B}(j+1+\alpha,\beta+n),\quad \alpha>-1,\, \beta>0,\\
\xi_n\nu_n^{j}\mathcal{M}w_{\rm CL}(j+1)&=n^{j+1+\alpha}\mathrm{B}(j+1+\alpha,n-j-2+\beta-\alpha),\quad \alpha>-1,\, \beta-\alpha>1.
}
Using the bound, for $j\in[1,n-1]$,
\Eq{}{
j^{j+1+\alpha}\mathrm{B}(j+1+\alpha,\beta+j+1)&\leq n^{j+1+\alpha}\mathrm{B}(j+1+\alpha,\beta+n)
}
and the bound, for $j\in \left[\max\{2(\beta-2)-\alpha,0\}, n-1\right]$,
\Eq{}{
\Gamma(j+1+\alpha)&\leq n^{j+1+\alpha}\mathrm{B}(j+1+\alpha,n-j-2+\beta-\alpha),
}
one can choose, respectively, the function $h$ and $j_*$, to be
\Eq{}{
h(j)&=h_{\rm Lag}(j):=\Gamma(j+1+\alpha),\quad j_*=\max\{1,\lceil1-\alpha\rceil\},\\
h(j)&=h_{\rm Jac}(j):=j^{j+1+\alpha}\mathrm{B}(j+1+\alpha,\beta+j+1),\quad j_*=\max\{1,\lceil1-\alpha\rceil,\lceil1+\beta\rceil\},\\
h(j)&=h_{\rm CL}(j):=\Gamma(j+1+\alpha),\quad j_*=\max\{1,\lceil1-\alpha\rceil,\lceil 2(\beta-2)-\alpha\rceil\}.
}
In each case, $h$ has the root asymptotic \eqref{eq: root asymp bound h}. Furthermore, with those particular choices of $j_*$, the functions $h_{\rm Lag}, h_{\rm Jac}, h_{\rm CL}$ are strictly increasing on $[j_*,n-1]$ and bounded from below by $1$.
}

\rem{rem: closure}{
It is interesting to notice that the set of P\'olya ensembles verifying  \autoref{assum:exist HE} enjoys the same closure property as the set of all P\'olya ensembles. Indeed, given two complex square matrices $X$ and $Y$ drawn from different (or same) Pólya ensembles whose respective P\'olya weight $w_X$ and $w_Y$ verify \autoref{assum:exist HE}, then, the matrix $XY$ also follows a Pólya ensemble, whose corresponding Pólya weight is given by the multiplicative convolution
\Eq{}{
w_{XY}(x):=(w_X * w_Y)(x)=\int_0^\infty w_X\left(\frac{x}{y}\right)w_Y(y)\frac{dy}{y}.
}
Due to the properties of the Mellin transform, we have
\Eq{}{
\mathcal{M}w_{XY}(s)= \mathcal{M}w_X(s)\, \mathcal{M} w_Y(s),
}
and it is then easy to see that the P\'olya weight $w_{XY}$ also verifies \autoref{assum:exist HE} with the scalings
\Eq{}{
\nu_{n,XY}=\nu_{n,X}\,\nu_{n,Y},\quad \xi_{n,XY}=\xi_{n,X}\,\xi_{n,Y},
}
where $\nu_{n,X},\xi_{n,X}$ (resp. $\nu_{n,Y},\xi_{n,Y}$) are the scalings associated with the P\'olya ensemble of $X$ (resp. Y). We note that the existence of the first point-wise limit in \eqref{eq: def limit polya w} is actually implied for almost all  $x\in \mathbb{R}_+$ by the other assumptions. However, we do not show it here.
}
Comparing \autoref{assum:exist HE PE} and \autoref{assum:exist HE},the first assumption \eqref{eq: def limit polya w} implies \eqref{eq: def limit kernel} through \autoref{prop: kernel HE polya}. The second assumption \eqref{eq: bound Tildew} implies \eqref{eq: bound p_j}. Indeed, using the explicit expression of the polynomials $p_j$ \eqref{eq:polya ortho functions} (rescaled by $1/\xi_n$), and the binomial identity, one can get the following bound, for $j=0,\ldots,n-1$,
\Eq{eq: bound pj polya}{
\abs{\frac{1}{\xi_n}p_j\left(\frac{y}{n\nu_n} \right) }=\abs{\sum_{c=0}^{j} \binom{j}{c} \frac{1}{\xi_n\nu_n^c \mathcal{M}w_n(c+1)}\left(\frac{-y}{n}\right)^c}\leq C \left(1+\frac{\abs{y}}{n} \right)^j\leq C \exp(\abs{y}),
}
with $C$ some $n$-independent positive constant. Assumption \eqref{eq: dom Tildew} is likely to imply \eqref{eq: bound q_j}, yet it is neither clear nor shown. Regarding the last assumption \eqref{eq: lower bound h} and \eqref{eq: root asymp bound h}, we show in \autoref{Appendix} that it implies \eqref{eq: asymp pj polynom}. This last assumption could possibly be lifted due to some property of the Pólya frequency functions such as their log-concavity \cite{Schoenberg1951,Foerster2020,Groechenig2020}. However, due to assumption \eqref{eq: dom Tildew}, it is not entirely clear if the set of Pólya ensembles verifying \autoref{assum:exist HE} is strictly included in the set of polynomial ensembles verifying \autoref{assum:exist HE PE}. Notwithstanding, in both cases, the formulas \eqref{eq:1,kpt poly inf}, \eqref{eq:1pointreal SV} and \eqref{eq:1pointreal} hold.  Hence the following theorem.

\theo{theo: 1,k polya ensemble}{Let $k\in \mathbb{N}$ and consider a random matrix drawn from a bi-unitarily invariant ensemble on $\mathrm{GL}(n,\mathbb{C})$  having a Pólya ensemble with joint probability density~\eqref{eq: DPP poly ens} for the squared singular values. Under \autoref{assum:exist HE}, instead of \autoref{assum:exist HE PE}, the results of \autoref{theo: 1,k poly ensemble} hold.
}
While the formulas in \autoref{theo: 1,k poly ensemble} seem to be quite general and are expressed explicitly in terms of the limiting kernel, they involve a double non-compact integral, which might not be easily computed. Yet, the structure of Pólya ensembles enables to go further by computing this double integral.

Let us first recall that, for Pólya ensembles, the kernel admits an integral representation \cite{Kuijlaars2014,Kieburg2016},
\Eq{eq: integ kernel polya}{
K_n(x,y)=\sum_{j=0}^{n-1} q_j(x) p_j(y)=n\int_0^1 dt\, q_n(xt)p_{n-1}(yt),
}
with the bi-orthonormal sets of functions $\left\{q_{j}\right\}_{j=0}^{n-1}$ and $\left\{p_{j}\right\}_{j=0}^{n-1}$  given by
\Eq{eq:polya ortho functions}{
p_{j}(x) &=\sum_{c=0}^{j} \binom{j}{c} \frac{(-x)^c}{\mathcal{M}w_n(c+1)}, \quad q_j(x) = \frac{1}{j!}\partial_x^j[x^j w_n(x)].}
For this result, the $n$-times differentiability of $w_n$ is needed, i.e., $w_n\in C^n(\mathbb{R})$, to have a well-defined function $q_n$. If the weight is only $(n-1)$-differentiable one needs to modify the formula to
\begin{equation}
K_n(x,y)=n\int_0^1 dt\, q_{n-1}(xt)p_{n-2}(yt)+q_{n-1}(x)p_{n-1}(y).
\end{equation}
One can rescale $q_n$ and $p_{n-1}$ as follows, where $q_n$ is expressed with the help of an inverse Mellin transform (the regularization function for the complex integral is omitted here),
\Eq{eq:polya ortho scaled1}{
Q_n(x) := \frac{\xi_n}{ \nu_n}q_n\left(\frac{x}{n \nu_n}\right)= \int_{\mathcal{C}_1} \frac{ds}{2 i\pi}\frac{\Gamma(n+1-s)n^{s-1} }{\Gamma(n)\Gamma(1-s)}\xi_n\nu_n^{s-1}\mathcal{M}w_n(s)x^{-s},
}
with $\mathcal{C}_1:=1+i\mathbb{R}$, and
\Eq{eq:polya ortho scaled2}{
P_{n-1}(y) :=\frac{1}{\xi_n}p_{n-1}\left(\frac{y}{n \nu_n}\right)=\sum_{c=0}^{n-1} \binom{n-1}{c} \frac{(-y)^c}{n^c \xi_n\nu_n^c\mathcal{M}w_n(c+1)}.
}
The scaled version of the kernel \eqref{eq: integ kernel polya} is then 
\Eq{eq: integ scaled kernel polya}{
 \frac{1}{n \nu_n}K_n\left(\frac{x}{n \nu_n},\frac{y}{n \nu_n}\right)=\int_0^1 dt Q_n(xt) P_{n-1}(yt),
}
for which the scaling limit is given by the following proposition, which is a modified version of \cite[Prop.2.8]{Kieburg2022}, with the corresponding \autoref{assum:exist HE} accounting for the additional scaling $\nu_n$. We relax the constraint on the domain of $y$ as the proposition is true for any $y\in\mathbb{R}$.

%---prop: kernel HE polya--------------------------------------------------------------------------
\prop{prop: kernel HE polya}{(\cite[Prop.2.8 modified]{Kieburg2022})
The hard edge limit of the kernel \eqref{eq: integ kernel polya} for Pólya ensembles satisfying the \autoref{assum:exist HE} is given by the point-wise limit
\Eq{}{
 K^{(\infty)}(x,y)&:=\lim_{n\to \infty} \frac{1}{n \nu_n}K_n\left(\frac{x}{n \nu_n},\frac{y}{n \nu_n}\right)=\int_0^1 dt Q^{(\infty)}(xt) P^{(\infty)}(yt), 
}
for any fixed $x\in\mathbb{R}_+$ and $y\in\mathbb{R}$, where, 
\Eq{eq:P inf}{
P^{(\infty)}(y):=\lim_{n\to\infty}P_{n-1}(y)=\sum_{j=0}^\infty \frac{(-y)^j}{j!\, \Tilde{w}^{(\infty)}(j+1)}.
}
\Eq{eq:Q inf}{
Q^{(\infty)}(x):=\lim_{n\to\infty}Q_n(x)=\int_{\mathcal{C}_1(\theta)} \frac{ds}{2 i\pi}\frac{\Tilde{w}^{(\infty)}(s)}{\Gamma(1-s)}x^{-s}
}
with $\mathcal{C}_1(\theta)=\left\{ 1+ie^{\sign(\tau)i \theta}\tau,\, \tau\in \mathbb{R} \right\}$ and $\theta \in ]0,\frac{\pi}{2}[$.

}
\rem{}{ The definition of the contour $\mathcal{C}_1(\theta)$ plays the role of a regularisation function, ensuring the convergence of the integral independently of $w^{(\infty)}$ \cite[Proof of Lemma 2.5]{Kieburg2022}. Note that while the definition of the $Q_n$ in \cite[Prop.2.8]{Kieburg2022} do not include the factor $\Gamma(n+1-s)n^{s-1}/\Gamma(n)$, they share the same scaling limit $Q^{(\infty)}$, by the same arguments used in the proof of \cite[Prop.2.8]{Kieburg2022}. 
}

It is known that for classical Pólya ensembles like Laguerre and Jacobi ensembles, the limiting kernel at the hard edge at the origin can be expressed in terms of the Bessel kernel, which admits the following integral representation (cf. \cite[Eq. (2.2)]{Tracy1994a}).

 %-------prop:Bessel kernel-------------------------------------------------------------
\prop{prop:Bessel kernel}{For all $x,y\in \mathbb{R}$,
\Eq{}{
 K_{\alpha}^{\mathrm{Bes}}(x,y):=&\frac{1}{4}\int_0^1 dt J_{\alpha}(\sqrt{xt}) J_{\alpha}(\sqrt{yt})=\frac{\sqrt{y} J'_{\alpha}(\sqrt{y})J_{\alpha}(\sqrt{x})-\sqrt{x} J'_{\alpha}(\sqrt{x})J_{\alpha}(\sqrt{y}) }{2(x-y)},
}
with $J_{\alpha}$ the usual Bessel function of the first kind of order $\alpha$.
}
\rem{}{In the formula \eqref{eq:1,kpt poly inf}, the second argument of the limiting kernel is negative. This is not an issue. Albeit $J_{\alpha}$ is analytic on $\mathbb{C}\setminus]-\infty,0]$, for negative arguments, one can use the relation with the modified Bessel function $I_\alpha$
\Eq{}{
\forall x \in \mathbb{R}_+,\quad J_\alpha(ix)=i^\alpha I_\alpha(x),\quad J'_\alpha(ix)=i^{\alpha-1}I'_\alpha(x).
}
For instance for $x,y \in \mathbb{R}_+$,
\Eq{}{
 K_{\alpha}^{\mathrm{Bes}}(x,-y)=i^\alpha \frac{\sqrt{y} I'_{\alpha}(\sqrt{y})J_{\alpha}(\sqrt{x})-\sqrt{x} J'_{\alpha}(\sqrt{x})I_{\alpha}(\sqrt{y}) }{2(x+y)}
 }
 and
\Eq{}{K_{\alpha}^{\mathrm{Bes}}(-x,-y)=(-1)^\alpha \frac{\sqrt{y} I'_{\alpha}(\sqrt{y})I_{\alpha}(\sqrt{x})-\sqrt{x} I'_{\alpha}(\sqrt{x})I_{\alpha}(\sqrt{y}) }{2(y-x)}.
}
}
In general, the Bessel kernel appears whenever the scaled Pólya weight converges to the Laguerre weight in the double scaling limit. 

%-------prop: limit kernel classic Pólya-------------------------------------------------------------
\prop{prop: limit kernel classic Pólya}{
If $\Tilde{w}^{(\infty)}(s)=\Gamma(s+\alpha)$, with $\Re(\alpha)>-1$, then, 
\Eq{}{
 K^{(\infty)}\left(x,y\right)&=\left(\frac{x}{y}\right)^{\frac{\alpha}{2}} \int_0^1 dt J_{\alpha}(2\sqrt{xt}) J_{\alpha}(2\sqrt{yt})=4 \left(\frac{x}{y}\right)^{\frac{\alpha}{2}}K_{\alpha}^{\mathrm{Bes}}(4x,4y)
}
with $K^{(\infty)}$ defined as in \autoref{prop: kernel HE polya}.
}
\begin{proof}[Proof of \autoref{prop: limit kernel classic Pólya}]
Using the definition of $P^{(\infty)}$ and $Q^{(\infty)}$ given in \autoref{prop: kernel HE polya} we get
\Eq{eq:P inf Gamma}{
P^{(\infty)}(y)=\sum_{j=0}^\infty \frac{(-y)^j}{j!\, \Gamma(j+1+\alpha)}=y^{-\frac{\alpha}{2}}J_{\alpha}(2\sqrt{y})
}
and 
\Eq{}{
Q^{(\infty)}(x)=\int_{\mathcal{C}_1(\theta)} \frac{ds}{2 i\pi}\frac{\Gamma(s+\alpha)}{\Gamma(1-s)}x^{-s}.
}
In this case the integral converges even for $\theta=0$, yielding $\mathcal{C}_1(0)=1+i\mathbb{R}=\mathcal{C}_1$. Thus, we get
\Eq{eq:Q inf Gamma}{
Q^{(\infty)}(x)=\int_{\mathcal{C}_1} \frac{ds}{2 i\pi}\frac{\Gamma(s+\alpha)}{\Gamma(1-s)}x^{-s}=x^{\frac{\alpha}{2}}J_{\alpha}(2\sqrt{x}).
}
By \autoref{prop:Bessel kernel} the result follows.
\end{proof}

Using \autoref{prop: kernel HE polya}, the scaling limit of the $1,k$-cross-covariance density function, is then given by the following corollary, which can be seen as the double scaling limit version of \autoref{prop:poly ensemble explicit}.
%=============CORO POLYA ENSEMBLE==============
\coro{coro:polya ensemble}{
For Pólya ensembles verifying the assumptions \autoref{assum:exist HE}, the scaling limit of the $1,k$-cross-covariance density functions between one squared eigenradius and $k$ squared singular values is given by \eqref{cov.prop inf} where $\mathbf{C}^{(\infty)}$ can be recast into the form
\Eq{cov.prop2 inf}{
\mathbf{C}^{(\infty)}(r;a_1,a_2)=& \sum_{\gamma=0,1} H_\gamma^{(\infty)}(r,a_1)\left[r^{-1}e^{-\frac{a_2}{r}}\left(\frac{a_2}{r}-1\right)^{1-\gamma}- V_\gamma^{(\infty)}(r,a_2) \right]\\
=&\, \partial_r\left[H_0^{(\infty)}(r,a_1)\left(e^{-\frac{a_2}{r}}-r V_1^{(\infty)}(r,a_2)\right) \right],\\
}
where, for $\gamma=0,1$,
\begin{eqnarray}
H_\gamma^{(\infty)}(x,y)&:=&\int_0^1 du\ Q^{(\infty)}(yu) \partial_{u}^\gamma \left[ u^\gamma\label{Hinf.def} \frac{\rho_{\rm EV}^{(\infty)}(xu)}{w^{(\infty)}(xu)}\right] , \\
V_\gamma^{(\infty)}(x,y) &:=&  \int_0^1 du\   P^{(\infty)}\left(yu\right) \left(u \partial_{u}\right)^{1-\gamma} w^{(\infty)}(xu),\label{Vinf.def}
\end{eqnarray}
with $P^{(\infty)}$ and $Q^{(\infty)}$ given respectively by \eqref{eq:P inf} and \eqref{eq:Q inf}, and 
\Eq{eq: rho EV inf polya}{
  \rho_{EV}^{(\infty)}(x)&= w^{(\infty)}(x)\sum_{k=0}^\infty \frac{x^k}{\Tilde{w}^{(\infty)}(k+1)}.
}
}

\begin{proof}[Proof of  \autoref{coro:polya ensemble}]
Let us first show that the scaling limit $\rho_{EV}^{(\infty)}$ is indeed given by \eqref{eq: rho EV inf polya}. From \cite[Eq. (4.7)]{Kieburg2016} the scaling version of the $1$-point function on the squared eigenradii is given by
\Eq{}{
\tfrac{n}{\nu_n}\rho_{EV}\left(\tfrac{x}{\nu_n} \right)=\frac{\xi_n}{\nu_n}w_n\left(\frac{x}{\nu_n} \right)  \sum_{j=0}^{n-1} \frac{x^j}{\xi_n\nu_n^{j}\mathcal{M}w_n(j+1)}.
}
Using the bound assumption \eqref{eq: lower bound h}, one gets the following bound
\Eq{eq: bound rho EV}{
  \sum_{j=0}^{n-1} \frac{x^j}{\xi_n\nu_n^{j}\mathcal{M}w_n(j+1)}\leq \sum_{j=0}^{j_*-1} \frac{x^j}{\xi_n\nu_n^{j}\mathcal{M}w_n(j+1)}+\sum_{j=j_*}^{n-1} \frac{x^j}{h(j)},
}
where we remind that $j_*$ is $n$-independent. By Cauchy-Hadamard theorem, due to the root asymptotic \eqref{eq: root asymp bound h}, the second sum on the RHS is a partial sum of a convergent series with infinite radius of convergence. The first $j_*$ summands are bounded by a constant independent of $n$ by assumption \eqref{eq: bound Tildew}. Moreover, as $h$ is also $n$-independent, this provides an upper bound to apply dominated convergence theorem. Thus, one is allowed to take the limit $n\to\infty$ inside the sum on the LHS and using the limits \eqref{eq: def limit polya w} yields the claim.

For the remaining part of the claim, we start from \autoref{prop:poly ensemble explicit} and apply Jacobi's formula for derivative of determinant, to get
\Eq{}{
\mathrm{cov}_{1,k}(r;a_1,\ldots,a_k)=&\frac{(n-k)!}{(n-1)!}\sum_{l,m=1}^k(-1)^{l+m}\mathbf{C}_n(r;a_l,a_m)\det[K_n(a_b,a_c)]_{\substack{b\in\llbracket 1,k\rrbracket\setminus\{l\}\\c\in\llbracket 1,k\rrbracket\setminus\{m\}}}.
}
One then proceeds with the appropriate scalings. The integrands in \eqref{V.def} and \eqref{H.def} are bounded by $n$-independent integrable functions on $]0,1]$, as it can be seen with \eqref{eq: bound rho EV}, \eqref{eq: bound pj polya} and with the bound given in \cite[Corollary 2.6]{Kieburg2022}. By Lebesgue's dominated convergence theorem, the limit $n\to \infty$ can be taken inside the integrals. Using, \autoref{prop: kernel HE polya}, all limits exist and one arrives at the first expression in \eqref{cov.prop2 inf}.

Then, using the relations
\begin{eqnarray}
H_1^{(\infty)}(x,y)&=&\partial_x \left[x H_0^{(\infty)}(x,y)\right], \\
V_0^{(\infty)}(x,y) &=&  x\partial_x V_1^{(\infty)}(x,y),
\end{eqnarray}    
which are the scaled versions of the relation given in \cite[Sec. 4.1]{Allard2025b}, finishes the proof.
\end{proof}
\ex{}{
When $\Tilde{w}^{(\infty)}(s)=\Gamma(s+\alpha)$ (cf. \autoref{ex:classic polya}), with $\Re(\alpha)>-1$, $P^{(\infty)}$ and $Q^{(\infty)}$ reduce, respectively, to \eqref{eq:P inf Gamma} and \eqref{eq:Q inf Gamma} and
\begin{eqnarray}
\rho_{EV}^{(\infty)}(x)&=& \frac{\gamma(\alpha,x)}{\Gamma(\alpha)},\\
H_0^{(\infty)}(x,y)&=&x^{-\alpha}y^{\alpha/2}\int_0^1 du\,J_{\alpha}(2\sqrt{y u})\, u^{-\alpha/2}  e^{x u} \frac{\gamma(\alpha,xu)}{\Gamma(\alpha)}, \\
V_1^{(\infty)}(x,y) &=&  x^{\alpha}y^{-\alpha/2}\int_0^1 du\,J_{\alpha}(2\sqrt{y u})\, u^{\alpha/2}  e^{-x u}, 
\end{eqnarray}    
where $(\alpha,x)\mapsto\gamma(\alpha,x)$ is the lower incomplete gamma function. Note that there is no issue at $\alpha=0$ as the function $(\alpha,z)\mapsto \frac{\gamma(\alpha,z)}{\Gamma(\alpha)}$ is holomorphic on $\mathbb{C}^2$.
}
%---------------------------------------------------------------------------------------
\subsection{Example: Muttalib-Borodin ensembles}
The Laguerre type Muttalib-Borodin ensemble with parameter $\theta'>0$ and weight function $x^{\alpha'}e^{-x}$, $\alpha'>-1$ is the following probability density function on $\mathbb{R}_+^n$ \cite{Muttalib1995,Borodin1998}
\Eq{}{
\frac{1}{Z_n}\Delta_n(x)\Delta_n(x^{\theta'})\det(x)^{\alpha'}e^{-\Tr(x)}.
}
Under the change of variables
\Eq{}{
x\mapsto x^{1/\theta'}, \quad \theta=\frac{1}{\theta'}, \quad \alpha=\frac{\alpha'+1}{\theta'}-1,
}
it becomes a Pólya ensemble with weight $w_{\rm MB}(x)=x^\alpha e^{-x^\theta}$, where $\theta>0$ and $\alpha>-1$. As the weight is independent of $n$, it remains the same in the large $n$ limit. Its Mellin transform is given by
\Eq{}{
\mathcal{M}w_{\rm MB}(s)=\Tilde{w}_{\rm MB}^{(\infty)}(s)=\frac{1}{\theta }\,\Gamma \left(\frac{\alpha+s}{\theta }\right).
}
Therefore, from \eqref{eq:biorthogonal}, the functions $p_{n-1}$ and $q_n$ composing the kernel $K_n$ \eqref{eq: integ kernel polya} are given by
\Eq{}{
p_{n-1}(y)=\theta\sum_{c=0}^{n-1} \binom{n-1}{c} \frac{(-y)^c}{\Gamma \left(\frac{\alpha+c+1}{\theta }\right)}
}
and 
\Eq{}{
q_{n}(x)=\frac{w_{\rm MB}(x)}{n!}\sum_{c=0}^{n} \frac{(x^\theta)^c}{c!}\sum_{l=0}^c \binom{c}{l}(-1)^l \frac{\Gamma\left(\theta l+\alpha+1+n\right)}{\Gamma\left(\theta l+\alpha+1\right)}.
}
By properties of the $\Gamma$ function, one can easily check that first assumptions in \autoref{assum:exist HE} are verified. For the last assumption, in this case $\nu_n=\xi_n=1$, and one can take
\Eq{}{
h(j)=h_{\rm MB}(j):=\frac{1}{\theta }\,\Gamma \left(\frac{j+1+\alpha}{\theta }\right),\quad j_*=\max\{\lceil2\theta-1\rceil,\lceil2\theta-1-\alpha\rceil\}.
}
The respective scaling limits $P^{(\infty)}$ and $Q^{(\infty)}$ (cf. \autoref{prop: kernel HE polya}) are given in terms of the Wright's generalization of the Bessel function
\Eq{eq: Wright B}{
J_{a,b}(x)=\sum_{j=0}^\infty \frac{(-x)^j}{j!\, \Gamma(aj+b)}.
}
Explicitly, 
\Eq{}{
P^{(\infty)}(y)=\theta J_{\frac{1}{\theta},\frac{\alpha+1}{\theta}}(y),\quad Q^{(\infty)}(x)=x^\alpha J_{\theta,\alpha+1}(x^\theta),\\
}
From \autoref{prop: kernel HE polya} one recovers the known limiting kernel
\Eq{}{
K^{(\infty)}_{\rm MB}(x,y)= \theta \int_0^1 dt (xt)^\alpha J_{\theta,\alpha+1}((xt)^\theta) J_{\frac{1}{\theta},\frac{\alpha+1}{\theta}}(yt).
}
One can then use \autoref{theo: 1,k polya ensemble}, by simply plugging the limiting kernel in \eqref{eq:1,kpt poly inf}, to obtain the $1,k$-point function, or, equivalently, in \eqref{cov.prop inf} to get the $1,k$-cross-covariance density.

Alternatively, one can use \autoref{coro:polya ensemble} to obtain the $1,k$-cross-covariance and consequently, the $1,k$-point function, with 
\begin{eqnarray}
\rho_{EV}^{(\infty)}(x)&=&\theta x^\alpha e^{-x^\theta} E_{\frac{1}{\theta},\frac{\alpha+1}{\theta}}(x) ,\\
H_0^{(\infty)}(x,y)&=&\theta \int_0^1 du\,(yu)^\alpha J_{\theta,\alpha+1}\left((yu)^\theta\right) E_{\frac{1}{\theta},\frac{\alpha+1}{\theta}}(xu), \\
V_1^{(\infty)}(x,y) &=& \theta\int_0^1 du\,J_{\frac{1}{\theta},\frac{\alpha+1}{\theta}}(yu)\, (xu)^{\alpha}  e^{-(x u)^\theta}, 
\end{eqnarray} 
where $E_{a,b}$ is the Mittag-Leffler function
\Eq{eq: Mittag}{
E_{a,b}(z):=\sum _{j=0}^{\infty }\frac{z^{j}}{\Gamma (a j+b )},\quad a>0.
}
This second formulation of the $1,k$-point function (and the $1,k$-cross-covariance) has the advantage of involving a single compact integral and is numerically much more efficient to use. One can even go further and use the explicit expressions of the Wright generalized Bessel function \eqref{eq: Wright B} and the Mittag-Leffler function \eqref{eq: Mittag} to compute the integrals $H_0^{(\infty)}(x,y)$ and $V_1^{(\infty)}(x,y)$ term by term.

Finally, one can check that \autoref{assum:exist HE PE} are verified as well when $\alpha\geq0$. Explicitly, using the asymptotic expansion of the Wright's function \cite{Wright1935}, the functions $f$ and $g$ can be taken to be
\Eq{}{
g(y)=\abs{y}+c_0,\qquad f(x)=c_1 x^{\frac{1}{1+\theta}}+c_2,\quad c_0,c_1,c_2>0,
}
and the last assumption is checked in \autoref{Appendix}.

%-----------------------------------------------------------
\section{Soft-hard edge scaling limit for Jacobi ensembles}\label{Soft-hard edge scaling limit for Jacobi ensembles}

While the density of most polynomial ensembles has a scaling limit at the origin, they often have only one hard edge on the singular values and eigenvalues. The Jacobi ensemble \eqref{eq:Jac buinv} with fixed parameter $\alpha$ and $\beta$, and corresponding P\'olya weight given in \eqref{eq: polya weights}, has a second hard edge for the singular values around the upper edge of its support. Let us study the corresponding double scaling limit around this other edge.

Using the Christoffel-Darboux formula, the kernel can be written
\Eq{eq:CD formula}{K_n(x,y)=a_n x^{\alpha}(1-x)^{\beta}\frac{P_{n}^{(\alpha,\beta)}(1-2x)P_{n-1}^{(\alpha,\beta)}(1-2y)-P_{n-1}^{(\alpha,\beta)}(1-2x)P_{n}^{(\alpha,\beta)}(1-2y)}{x-y}
}
with the usual Jacobi polynomials
\Eq{}{
P_{n}^{{(\alpha ,\beta )}}(z)=\frac{\Gamma(\alpha +n+1)}{n!\,\Gamma(\alpha +\beta +n+1)}\sum _{m=0}^{n}\binom{n}{m}\frac{\Gamma(\alpha +\beta +n+m+1)}{\Gamma(\alpha +m+1)}\left({\frac  {z-1}{2}}\right)^m 
}
and the coefficient
\Eq{eq: coeff an}{
a_n=-\frac{ \Gamma(n+1) \Gamma(n+\alpha +\beta +1)}{(\alpha +\beta +2 n) \Gamma(n+\alpha ) \Gamma(n+\beta )}.
}
\rem{}{Note that the Jacobi polynomials associated with our Jacobi ensemble are orthogonal on $[0,1]$ and are therefore the polynomials in \eqref{eq:CD formula}. We express them in terms of the usual Jacobi polynomials to use known asymptotics.
}
One can easily show, using the reflection property of the Jacobi polynomials
\Eq{}{
P_{n}^{(\alpha ,\beta )}(-z)=(-1)^{n}P_{n}^{(\beta ,\alpha )}(z)
}
and \autoref{prop: kernel HE polya}, that the limiting kernel for the point process on the singular values is 
\Eq{eq: limiting kernel Jac UHE}{
 \lim_{n\to \infty} \frac{1}{n^2}K_n\left(1-\frac{x}{n^2},1-\frac{y}{n^2}\right)&=4 \left(\frac{x}{y}\right)^{\frac{\beta}{2}}K_{\beta}^{\mathrm{Bes}}(4x,4y),
}
therefore, the limiting microscopic level density on the singular values is
\Eq{}{
 \lim_{n\to \infty} \rho_{\rm SV}\left(1-\frac{a}{n^2}\right)=J_{\beta }\left(2 \sqrt{a}\right)^2-J_{\beta -1}\left(2 \sqrt{a}\right) J_{\beta +1}\left(2 \sqrt{a}\right).
}
However, the limit of the level density of the eigenradii around $r=1$ is trivially null as the limiting macroscopic level density on the squared eigenradii is given by \cite[Eq.(3.12)]{Akemann2014}
\Eq{eq: limit rhoEV}{
\lim_{n\to \infty} \rho_{\rm EV}(r)=\frac{1}{(1-r)^2} \Theta \left(\frac{1}{2}-r\right).
}
Note that the Heaviside function in \eqref{eq: limit rhoEV} appears only in the limit and is not associated with a hard edge. Actually, if one zooms on a scale $1/\sqrt{n}$ around $r=1/2$, the limiting microscopic level density of the eigenradii is given in terms of a complementary error function \cite[Eq.(3.36)]{Akemann2014}. There is, therefore, no other hard edge for the eigenradii than the one at the origin. The upper edge around $1/2$ is a soft edge. This is due to the $n$ dependence of the P\'olya weight, namely, the repulsion from the upper edge grows with $n$, in this case, for the eigenvalues, cf. \eqref{eq: polya weights}.

The $1,k$-cross-covariance density then admits a scaling limit around the soft edge of the eigenradii and the upper hard edge of the singular values given in the following theorem. We will refer to this double scaling limit as Soft-Hard (SH) edge limit.
%=============THEOREM JAC UP HE==============
\theo{theo:jac up HE}{
For the Jacobi ensemble with P\'olya weight $w_{\rm Jac}(x)=x^\alpha (1-x)^{\beta+n-1} \Theta(1-x)$, with $\alpha,\beta$ fixed and $\alpha>-1$, $\beta>0$, the point-wise limit 
\Eq{}{
\mathrm{cov}_{1,k\, {\rm SH}}^{(\infty)}(r; a_1,\ldots,a_k):=\lim_{n\to\infty} n^{3/2-k}\,\mathrm{cov}_{1,k}\left(\frac{1}{2}-\frac{r}{\sqrt{n}}; 1-\frac{a_1}{n^2},\ldots,1-\frac{a_k}{n^2} \right)
}of the $1,k$-cross-covariance density functions between one squared eigenradius and $k$ squared singular values is given, for all fixed $(r,a_1,\ldots,a_k)\in \mathbb{R}\times \mathbb{R}_+^k$, by
\Eq{}{
&\mathrm{cov}_{1,k\, {\rm SH}}^{(\infty)}(r; a_1,\ldots,a_k)\\
&\quad = \frac{4}{\sqrt{\pi }} e^{-4 r^2} \partial_\mu \det\left[ \int_0^1 dt J_{\beta}(2\sqrt{a_b t}) J_{\beta}(2\sqrt{a_c t})+\mu J_{\beta}(2\sqrt{a_b }) J_{\beta}(2\sqrt{a_c }) \right]_{b,c=1}^k \Bigg|_{\mu=0}\\
&\quad =-\frac{4^k}{\sqrt{\pi }} e^{-4 r^2}\det\left(\begin{array}{c c} 
    	0\quad &  J_{\beta}(2\sqrt{a_c })  \\
     J_{\beta}(2\sqrt{a_b }) \quad & K_{\beta}^{\mathrm{Bes}}(4a_b,4a_c) 
\end{array}\right)_{b,c=1}^k.\\
}
In particular, in the case $k=1$, the limiting cross-covariance reads
\Eq{}{
\mathrm{cov}^{(\infty)}_{\rm SH}(r; a)= \frac{4}{\sqrt{\pi }} e^{-4 r^2} J_{\beta}(2\sqrt{a})^2.
}
}
\rem{rem: speed decrease 2}{
Echoing \autoref{rem: speed decrease 1}, the above result implies
\Eq{}{
\frac{1}{n^{2k+1/2}}\,\mathrm{cov}_{1,k}\left(\frac{1}{2}-\frac{r}{\sqrt{n}}; 1-\frac{a_1}{n^2},\ldots,1-\frac{a_k}{n^2} \right)\underset{n\to \infty}{=}\frac{1}{n^{k+2}}\left[\mathrm{cov}_{1,k\, {\rm SH}}^{(\infty)}(r;a_1,\ldots,a_k)+o(1)\right],
}
which means that the scaling limit of the $1,k$-cross-covariance, for Jacobi ensembles, at the soft-hard edge is of order $O(1/n^{k+2})$ as $n\to \infty$. It is therefore a factor $1/n$ smaller than the limiting $1,k$-point function and can be seen as the first correction term.
}

% %-----------------------------------------------------------
% \section{Proofs}\label{Proofs}
%-----------------------------------------------------------
\section{Proof: Hard edge at the origin}\label{Proof HE}
%-----------------------------------------------------------------------------
\subsection{Proof of \autoref{theo: 1,k poly ensemble}}
 The proof of \autoref{theo: 1,k poly ensemble} relies on Lebesgue's dominated convergence theorem, hence \autoref{assum:exist HE PE}. However, we do not assume the existence of an integrable upper bound, simply of a continuous one. One has, therefore, to split the integration into two parts: a compact part where one can use bounded convergence theorem to bring the limit below the integral and a non-compact part which vanishes in the appropriate limits. The most important step of the proof relies therefore on the following Lemma.
\lem{lem: remainder}{Let $k\in\mathbb{N}$. Under \autoref{assum:exist HE PE}, for almost all $r\geq 0$,
 \Eq{eq: lim remainder}{
&\lim_{R\to\infty}\lim_{n\to\infty}\left(\int_{0}^{R} dt \int_{C_R}^{n}dv+\int_{R}^{\infty} dt \int_0^{n}dv\right)\frac{n^2}{v} \varphi_n\left(\frac{v}{n }, \frac{t}{n}\right)\\
&\quad \ \times \det\left(\begin{array}{c c} 
    	 \frac{1}{n\nu_n}K_n\left(\frac{rv}{n\nu_n},\frac{-rt}{n\nu_n}\right)\quad & \frac{1}{n\nu_n}K_n\left(\frac{rv}{n\nu_n},\frac{a_c}{n\nu_n}\right)-\delta(rv-a_c)  \\
     \frac{1}{n\nu_n} K_n\left(\frac{a_b}{n\nu_n},\frac{-rt}{n\nu_n}\right) \quad & \frac{1}{n\nu_n}K_n\left(\frac{a_b}{n\nu_n},\frac{a_c}{n\nu_n}\right) 
\end{array}\right)_{b,c=1}^k=0,
}
with $C_R$ a positive constant depending only on $R$ and such that $\lim_{R\to \infty}C_R=\infty$.
}
Expanding the determinant in \eqref{eq: lim remainder} one can see there are 3 different types of integrands, each time integrated on 2 different domains. We therefore introduce the following notations for the 6 different quantities:
\Eq{eq: I0}{
I_{0,n}^-(R):=&\int_{R}^{\infty} dt \int_{0}^{n}\frac{dv}{v} n^2\varphi_n\left(\frac{v}{n }, \frac{t}{n}\right)\frac{1}{n\nu_n}K_n\left(\frac{ a_1}{n\nu_n},\frac{-rt}{n\nu_n}\right)\delta(rv-a_2),\\
I_{0,n}^+(R):=&\int_{0}^{R} dt \int_{C_R}^{n}\frac{dv}{v} n^2\varphi_n\left(\frac{v}{n }, \frac{t}{n}\right)\frac{1}{n\nu_n}K_n\left(\frac{ a_1}{n\nu_n},\frac{-rt}{n\nu_n}\right)\delta(rv-a_2),
}
\Eq{eq: I1}{
I_{1,n}^-(R):=&\int_{R}^{\infty} dt \int_{0}^{n}\frac{dv}{v} n^2\varphi_n\left(\frac{v}{n }, \frac{t}{n}\right)\frac{1}{n\nu_n}K_n\left(\frac{ r v}{n\nu_n},\frac{-rt}{n\nu_n}\right),\\
I_{1,n}^+(R):=&\int_{0}^{R} dt \int_{C_R}^{n}\frac{dv}{v} n^2\varphi_n\left(\frac{v}{n }, \frac{t}{n}\right)\frac{1}{n\nu_n}K_n\left(\frac{ r v}{n\nu_n},\frac{-rt}{n\nu_n}\right),
}
and
\Eq{eq: I2}{
I_{2,n}^-(R):=&\int_{R}^{\infty} dt \int_{0}^{n}\frac{dv}{v} n^2\varphi_n\left(\frac{v}{n }, \frac{t}{n}\right)\frac{1}{(n\nu_n)^2}K_n\left(\frac{rv}{n\nu_n},\frac{a_1}{n\nu_n}\right) K_n\left(\frac{a_2}{n\nu_n},\frac{-rt}{n\nu_n}\right),\\
I_{2,n}^+(R):=&\int_{0}^{R} dt \int_{C_R}^{n}\frac{dv}{v} n^2\varphi_n\left(\frac{v}{n }, \frac{t}{n}\right)\frac{1}{(n\nu_n)^2}K_n\left(\frac{rv}{n\nu_n},\frac{a_1}{n\nu_n}\right) K_n\left(\frac{a_2}{n\nu_n},\frac{-rt}{n\nu_n}\right).
}
The proof of \autoref{lem: remainder} relies on the more elementary lemma which follows.
\lem{lem: asymp qj}{Let $n\in \mathbb{N}$, $n>2$, and a continuous function $f:\mathbb{R}_+^*\to\mathbb{R}_+$ such that $f(u)=o(u)$, and $u\mapsto\exp(f(u))\in \mathrm{L}^1_{\rm loc}(\mathbb{R}_+)$. There exists $F:\mathbb{R}_+\to\mathbb{R}_+$ a $n$-independent,  continuous and strictly decreasing function such that $F(x)\to 0$, as $x\to +\infty$, and that we have
\Eq{}{
 \int_{x}^{n} dv (1+v) \left(1-\frac{v}{n}\right)^{n-2} \exp(f(v))\leq F(x),\quad \forall x\in[0,n].
 }
}
\begin{proof}
Using the fact 
\Eq{eq: bound exp}{
\forall v\in [0,n], \quad \left(1-\frac{v}{n} \right)^{n-2}
\leq 2e^{-v},
}
one has
\Eq{}{
\int_{x}^{n} dv\, (1+v) \left(1-\frac{v}{n}\right)^{n-2} \exp(f(v))&\leq 2\int_{x}^{n} dv\, (1+v)  \exp(-v+f(v))\\
&\leq 2\int_{x}^{\infty} dv\, (1+v)  \exp(-v+f(v)).
}
Due to the asymptotic behavior of $f$, $f(u)=o(u)$, as $u\to \infty$,
\Eq{}{
F(x):=2\int_{x}^{\infty} dv\, (1+v)  \exp(-v+f(v))
}
is the remainder of a converging integral of a continuous integrand on $\mathbb{R}_+^*$. Moreover, $F(0)$ is well defined due to $u\mapsto\exp(f(u))\in \mathrm{L}^1_{\rm loc}(\mathbb{R}_+)$. Therefore $F$ is continuous on $\mathbb{R}_+$ and $F(x)\to 0$, as $x\to +\infty$. Moreover, as the integrand is strictly positive, we have $F'<0$, which implies $F$ is strictly decreasing.
\end{proof}

\begin{proof}[Proof of \autoref{lem: remainder}]
Let us start with $I_{1,n}^-(R)$. The first step is to bound $\varphi_n$. Explicitly, we have
\Eq{}{
\frac{n^2}{v}\varphi_n\left(\frac{v}{n }, \frac{t}{n}\right)=\left(t+v-1-\frac{1}{n^2}v t\right)\left(1-\frac{v}{n }\right)^{n-2}\left(1+\frac{t}{n}\right)^{-(n+2)}.
}
Thus, we can choose the following bound
\Eq{}{
\forall(v,t)\in [0,n]\times \mathbb{R}_+,\quad \abs{\frac{n^2}{v}\varphi_n\left(\frac{v}{n }, \frac{t}{n}\right)}\leq (1+v) \left(1-\frac{v}{n }\right)^{n-2}\frac{(1+t)}{\left(1+\frac{t}{n}\right)^{n+2}}.
}
Regarding the kernel, using the expression \eqref{eq: Kernel poly ens} with the assumption  \eqref{eq: bound q_j} yields, for all $v$ and $t$,
\Eq{eq: bound kernel rem}{
\abs{\frac{1}{n\nu_n}K_n\left(\frac{ r v}{n\nu_n},\frac{-rt}{n\nu_n}\right)}\leq \frac{\exp(f(rv))}{n} \sum_{j=0}^{n-1} \abs{p_j\left(-\frac{r t}{n\nu_n} \right) },
}
where $f(u)=o(u)$, as $u\to \infty$, and $u\mapsto\exp(f(u))\in \mathrm{L}^1_{\rm loc}(\mathbb{R}_+)$. From there, we get
\Eq{}{
&\abs{I_{1,n}^-(R)}\\
&\leq  \int_{R}^{\infty} dt  \frac{ (1+t)}{\left(1+\frac{t}{n}\right)^{n+2}}\frac{1}{n}\sum_{j=0}^{n-1} \abs{p_j\left(-\frac{r t}{n\nu_n} \right) } \int_{0}^{n} dv (1+v) \left(1-\frac{v}{n}\right)^{n-2} \exp(f(rv)) .
}
As $r$ is fixed and $f(u)=o(u)$, as $u\to \infty$, and $\exp(f(u))=o(1/u)$, as $u\to 0$, using \autoref{lem: asymp qj} yields
\Eq{}{
\int_{0}^{n} dv (1+v) \left(1-\frac{v}{n}\right)^{n-2} \exp(f(rv))\leq F(0)=C,
}
for some positive constant $C$, independent of $n$. Finally, using assumption \eqref{eq: asymp pj polynom} finishes to show that $\lim_{R\to\infty}\lim_{n\to\infty}I_{1,n}^-(R)=0$.

Turning to $I_{1,n}^+(R)$, in a similar way, one has
\Eq{}{
&\abs{I_{1,n}^+(R)}\\
&\leq \int_{0}^{R} dt  \frac{ (1+t)}{\left(1+\frac{t}{n}\right)^{n+2}}\frac{1}{n}\sum_{j=0}^{n-1} \abs{p_j\left(-\frac{r t}{n\nu_n} \right) }  \int_{C_R}^{n} dv (1+v) \left(1-\frac{v}{n}\right)^{n-2} \exp(f(rv)).
}
Using \autoref{lem: asymp qj} and \eqref{eq: bound p_j} leads to
\Eq{eq: use C_R}{
\abs{I_{1,n}^+(R)}\leq F(C_R)\int_{0}^{R} dt  \frac{ (1+t)}{\left(1+\frac{t}{n}\right)^{n+2}}\exp\left(g(-rt)\right).
}
Recalling that $r$ is fixed, one can always find a constant $S_R$ independent of $n$ to get the very loose bound 
\Eq{}{
\int_{0}^{R} dt  \frac{ (1+t)}{\left(1+\frac{t}{n}\right)^{n+2}}\exp\left(g(-rt)\right)\leq R \sup_{t\in[0,R]}\exp\left(g(-rt)\right) \leq S_R
}
Then, as $F$ is strictly decreasing, it is invertible, and, for all $R$ large enough such that $\exp(-R-S_R)\in F(\mathbb{R_+})$, one can set 
\Eq{}{
C_R=F^{-1}\left(\exp(-R-S_R)\right).
}
The function $R\mapsto C_R$ is strictly increasing as it is the composition of two strictly decreasing functions. Moreover, since $ \exp(-R - S_R) \to 0 $ as $ R \to \infty $, and $ F(x) \to 0 $ as $ x \to \infty $, it follows from the continuity and monotonicity of $ F $ (by the continuous inverse theorem) that $ C_R \to \infty $ as $ R \to \infty $. This finishes to show that $\lim_{R\to\infty}\lim_{n\to\infty}I_{1,n}^+(R)=0$. Note, here, the need to choose the domains of the remainder integrals carefully \eqref{eq: domains}.

Moving to $I_{2,n}^-(R)$, $I_{2,n}^+(R)$, as $a_1, a_2$ are fixed, using once again the expression \eqref{eq: Kernel poly ens} with \eqref{eq: bound q_j} yields, for all $v$ in compact subsets of $\mathbb{R}_+$ and for all $t\in\mathbb{R}$,
\Eq{}{
&\abs{\frac{1}{(n\nu_n)^2}K_n\left(\frac{rv}{n\nu_n},\frac{a_1}{n\nu_n}\right) K_n\left(\frac{a_2}{n\nu_n},\frac{-rt}{n\nu_n}\right)}\\
&\leq \exp(g(a_1)+f(a_2)+f(rv))\frac{1}{n}\sum_{j=0}^{n-1} \abs{p_j\left(-\frac{r t}{n\nu_n} \right) }. 
}
This is the same upper bound as \eqref{eq: bound kernel rem} up to the multiplicative constant $\exp(g(a_1)+f(a_2))$ which is independent of $n$ and $R$. Thus, we also have $\lim_{R\to\infty}\lim_{n\to\infty}I_{2,n}^\pm(R)=0$.

Let us now consider $I_{0,n}^\pm(R)$. We immediately have $\lim_{R\to\infty}\lim_{n\to\infty}I_{0,n}^+(R)=0$, as a simple rewriting gives
\Eq{}{
I_{0,n}^+(R):=\Theta(n-a_2)\Theta(a_2-C_R)\int_{0}^{R} dt \frac{n^2}{a_2} \varphi_n\left(\frac{a_2}{rn }, \frac{t}{n}\right)\frac{1}{n\nu_n}K_n\left(\frac{ a_1}{n\nu_n},\frac{-rt}{n\nu_n}\right)
}
and, as $a_2$ is fixed, $I_{0,n}^+(R)$ vanishes for $C_R>a_2$. Regarding $I_{0,n}^-(R)$, one follows the same steps as for the case $I_{1,n}^-(R)$, using assumption \eqref{eq: asymp pj polynom},  this time without the $v$ integral, and gets $\lim_{R\to\infty}\lim_{n\to\infty}I_{0,n}^-(R)=0$.

Finally, to finish the proof, one expands the determinant in \eqref{eq: lim remainder} using Leibniz formula, uses the finite limit assumption \eqref{eq: def limit kernel} for the terms where none of the arguments are integrated and uses the vanishing limits of $I_{0,n}^\pm(R), I_{1,n}^\pm(R)$ and $I_{2,n}^\pm(R)$ computed above.
\end{proof}

\begin{proof}[Proof of \autoref{theo: 1,k poly ensemble}]
Starting from \eqref{eq:1,kpt poly}, proceeding with the scaling of the variables and doing the following change of variables $t\mapsto \frac{t}{n}$ and $v\mapsto \frac{rv}{n\nu_n}$ yields
\Eq{eq: scaled 1,k pt}{
&\frac{n^{k+1}}{\nu_n (n\nu_n)^{k}}f_{1,k}\left(\frac{r}{\nu_n}; \frac{a_1}{n\nu_n},\ldots, \frac{a_k}{n\nu_n}\right)\\
&\quad =\frac{n^k (n-k)!}{n!}\int_{0}^{\infty} dt \int_0^{n}\frac{dv}{v} n^2\varphi_n\left(\frac{v}{n }, \frac{t}{n}\right)\\
&\quad \ \times \det\left(\begin{array}{c c} 
    	 \frac{1}{n\nu_n}K_n\left(\frac{rv}{n\nu_n},\frac{-rt}{n\nu_n}\right)\quad & \frac{1}{n\nu_n}K_n\left(\frac{rv}{n\nu_n},\frac{a_c}{n\nu_n}\right)-\delta(rv-a_c)  \\
     \frac{1}{n\nu_n} K_n\left(\frac{a_b}{n\nu_n},\frac{-rt}{n\nu_n}\right) \quad & \frac{1}{n\nu_n}K_n\left(\frac{a_b}{n\nu_n},\frac{a_c}{n\nu_n}\right) 
\end{array}\right)_{b,c=1}^k,
}
where we have used the scaling property of the Dirac delta function
\Eq{}{
\forall c\neq 0, \quad \delta\left(\frac{x}{c}\right)=c \,\delta(x).
}
Let us now split the domain of integration $\mathcal{D}_t$ of the $t$ integral (resp. $\mathcal{D}_v$ for the $v$ integral) in two parts: one on which all the limits exist and one on which the integral will vanish as $n\to \infty$ then $R\to\infty$. The appropriate choice given \autoref{assum:exist HE PE} is the following
\Eq{eq: domains}{
\mathcal{D}_t=[0,R]\cup [R,\infty[,\quad \mathcal{D}_v=[0,C_R]\cup [C_R,n],
}
where $C_R$ is the constant in \autoref{lem: remainder} corresponding to the bound $f$ in \eqref{eq: bound q_j}. The choice of $C_R$ is only a technical details needed in the proof of \autoref{lem: remainder} (cf. \eqref{eq: use C_R}). Coming back to the integrand, we explicitly have
\Eq{}{
n^2\varphi_n\left(\frac{v}{n }, \frac{t}{n}\right)=v\left(t+v-1-\frac{1}{n^2}v t\right)\left(1-\frac{v}{n }\right)^{n-2}\left(1+\frac{t}{n}\right)^{-(n+2)}
}
and for any $(v,t)\in [0,C_R]\times [0,R]$, we have the uniform convergence of the limit
\Eq{}{
\lim_{n\to\infty} n^2\varphi_n\left(\frac{v}{n }, \frac{t}{n}\right)= v(t+v-1)e^{-v} e^{-t}=\varphi^{(\infty)}(v,t).
}
Using the pointwise convergence of the sequence of scaled kernels \eqref{eq: def limit kernel} and the bounds of the $p_j$ \eqref{eq: bound p_j} and $q_j$ \eqref{eq: bound q_j} one can use bounded convergence theorem to take the limit $n\to \infty$, below the integrals. Using the limit,
\Eq{eq: lim factorial}{
\lim_{n\to\infty} \frac{n^k (n-k)!}{n!}=1,
}
for fixed $k$ and then taking the limit $R\to\infty$, yields \eqref{eq:1,kpt poly inf}. Using \autoref{lem: remainder} finishes to prove the remainder vanishes in the large $n$, then $R$, limits.
\end{proof}

%-----------------------------------------------------------------------------
\subsection{Proof of \autoref{theo: 1,k polya ensemble}}
The idea of the proof of \autoref{theo: 1,k polya ensemble} is to use \autoref{coro:polya ensemble}, which is proven independently of \autoref{theo: 1,k polya ensemble}, and reintroduce the integrals in \eqref{cov.prop2 inf} to get back to the integral formulation \eqref{eq:1,kpt poly inf}. The proof of \autoref{theo: 1,k polya ensemble} can be made shorter by using the following lemma. Without it, one needs to proceed by integration by part.
\lem{lem: diff form}{Given $\varphi^{(\infty)}$ defined in \eqref{eq: def phi inf} and $K^{(\infty)}$ \eqref{eq: def limit kernel},
\Eq{eq: deriv identity1}{
\int_{0}^{\infty}dt  \int_0^\infty \frac{dv}{v} \varphi^{(\infty)}\left(v,t\right) K^{(\infty)}\left(rv,-r t\right)=\partial_r\int_{0}^{\infty}dt  \int_0^\infty dv\, e^{-v}e^{-t} r K^{(\infty)}\left(rv,-r t\right)
}
and
\Eq{eq: deriv identity2}{&\int _{0}^{\infty}dt  \int_0^\infty \frac{dv}{v} \varphi^{(\infty)}\left(\tfrac{v}{r},t\right) K^{(\infty)}\left(a_b,-rt\right)\left[\delta(v-a_c)-K^{(\infty)}(v,a_c)\right]\\
=&\,\partial_r \int _{0}^{\infty}dt  \int_0^\infty \frac{dv}{r}e^{-\frac{v+t}{r}}K^{(\infty)}\left(a_b,-t\right)\left[\delta(v-a_c)-K^{(\infty)}(v,a_c)\right].
}}
\begin{proof}
Rigorously, one needs to show that all integrals in \autoref{lem: diff form} are convergent. This can be shown using \autoref{prop: kernel HE polya}. The convergence is then given by the exponential decay of $\varphi^{(\infty)}$ \eqref{eq: def phi inf} and the last of \autoref{assum:exist HE}, which implies the asymptotic
\Eq{}{
P^{(\infty)}(\abs{x})=\exp(o(\abs{x})), \quad \abs{x}\to \infty,
}
as it can be seen in \eqref{eq: inequalities bound}. The explicit computation of those integrals is actually the object of the proof of \autoref{theo: 1,k polya ensemble}.

Starting with \eqref{eq: deriv identity1}, proceeding with a change of variables $t\mapsto t/r$ and $v\mapsto v/r$ and invoking Leibniz integral rule yields
\Eq{}{
&\partial_r\int_{0}^{\infty}dt  \int_0^\infty dv\, e^{-v}e^{-t} r K^{(\infty)}\left(rv,-r t\right)\\
=&\int_{0}^{\infty}dt  \int_0^\infty dv\, e^{-\frac{v+t}{r}} \left(\frac{t+v}{r^3}-\frac{1}{r^2} \right) K^{(\infty)}\left(v,-t\right).
}
Doing the inverse change of variables $t\mapsto rt$ and $v\mapsto rv$, on the RHS finishes to show \eqref{eq: deriv identity1}. The second equality \eqref{eq: deriv identity2} is done similarly.
\end{proof}
\begin{proof}[Proof of \autoref{theo: 1,k polya ensemble}]
As \autoref{coro:polya ensemble} is proven independently of \autoref{theo: 1,k polya ensemble}, we need to show that the RHS of \eqref{eq: Cinf integ} and  \eqref{cov.prop2 inf} are equal and that the RHS of \eqref{eq:1pointreal} and \eqref{eq: rho EV inf polya} are also equal. Let us first show that
\Eq{eq: id prove}{
  \rho_{EV}^{(\infty)}(r)&= w^{(\infty)}(r)\sum_{k=0}^\infty \frac{r^k}{\Tilde{w}^{(\infty)}(k+1)}=\int_{0}^{\infty}dt  \int_0^\infty \frac{dv}{v} \varphi^{(\infty)}\left(v,t\right) K^{(\infty)}\left(rv,-r t\right).
}
Using \autoref{lem: diff form}, along with \autoref{prop: kernel HE polya}, one can apply the derivative in $r$ using a change of variable $u\mapsto u/r$ to get
\Eq{eq: 1}{
&\int_{0}^{\infty}dt  \int_0^\infty \frac{dv}{v} \varphi^{(\infty)}\left(v,t\right) K^{(\infty)}\left(rv,-r t\right)\\
=&\,\partial_r\int_{0}^{\infty}dt  \int_0^\infty dv\, e^{-v}e^{-t} r \int_0^1 du\, Q^{(\infty)}(rvu) P^{(\infty)}(-rtu)\\
=&\int_{0}^{\infty}dt  \int_0^\infty dv\, e^{-v}e^{-t} Q^{(\infty)}(rv) P^{(\infty)}(-rt).
}
As the two integrals factorise, one can use the expression of $P^{(\infty)}$ \eqref{eq:P inf} to compute the $t$ integral, which is equal to $j!$. This simplifies with the factorial factor in the denominator and we get
\Eq{}{
 \int_{0}^{\infty}dt  \int_0^\infty dv\, e^{-v}e^{-t} Q^{(\infty)}(rv) P^{(\infty)}(-rt)=\int_0^\infty dv\, e^{-v} Q^{(\infty)}(rv)\sum_{j=0}^\infty \frac{r^j}{ \Tilde{w}^{(\infty)}(j+1)}.
}
Then, using the expression of $Q^{(\infty)}$ \eqref{eq:Q inf}, one would want to exchange the order of integration. However, it is not immediate. Instead we define,
\Eq{}{
W(r):=\int_0^\infty dv\, e^{-v} Q^{(\infty)}(rv).
}
The function $W$ is some multiplicative Mellin convolution between $G(x):=\exp(-x)$ and $Q^{(\infty)}$. Applying the Mellin transform on $W$ then yields
\Eq{}{
\mathcal{M}W(s)=\mathcal{M}G(1-s)\, \mathcal{M}Q^{(\infty)}(s),
}
where $\mathcal{M}G(1-s)=\Gamma(1-s)$ and
\Eq{}{
\mathcal{M}Q^{(\infty)}(s)= \frac{\mathcal{M}w^{(\infty)}(s)}{\Gamma(1-s)}.
}
Hence,
\Eq{}{
\mathcal{M}W(s)=\mathcal{M}w^{(\infty)}(s),
}
which yields, by the Mellin inversion theorem,
\Eq{}{
W(r)=w^{(\infty)}(r).
}
This finishes to show \eqref{eq: id prove}.

It remains to be shown that 
\Eq{}{
&\partial_r \int _{0}^{\infty}dt  \int_0^\infty \frac{dv}{r}e^{-\frac{v+t}{r}}K^{(\infty)}\left(a_b,-t\right)\left[\delta(v-a_c)-K^{(\infty)}(v,a_c)\right]\\
=&\,\partial_r\left[H_0^{(\infty)}(r,a_b)\left(e^{-\frac{a_c}{r}}-r V_1^{(\infty)}(r,a_c)\right) \right],
}
with $H_0^{(\infty)}$ and $V_1^{(\infty)}$ \eqref{Hinf.def} and \eqref{Vinf.def}. Following similar steps as the first part of the proof, one arrives at 
\Eq{}{
&\partial_r \int _{0}^{\infty}dt  \int_0^\infty \frac{dv}{r}e^{-\frac{v+t}{r}}K^{(\infty)}\left(a_b,-t\right)\left[\delta(v-a_c)-K^{(\infty)}(v,a_c)\right]\\
=&\,\partial_r\left[H_0^{(\infty)}(r,a_b)e^{-\frac{a_c}{r}}\right]-\partial_r\left[ r H_0^{(\infty)}(r,a_b) V_1^{(\infty)}(r,a_c) \right]
}
which is \eqref{cov.prop2 inf}. Finally, \eqref{eq:1pointreal SV} follows directly from \autoref{prop: kernel HE polya} and this finishes the proof.
\end{proof}
%-----------------------------------------------------------
\section{Proof: Soft-hard edge for Jacobi ensembles}\label{Proof SH}
%=============POLY ENSEMBLE==========================================================
\subsection{Proof of \autoref{theo:jac up HE}}
Let us start by defining
\Eq{}{
\mathbf{C}_n(r;a_1,a_2):=\int _{0}^{\infty}dt  \int_0^r \frac{dv}{v} \varphi_n\left(\tfrac{v}{r},t\right) K_n\left(a_1,-rt\right)\left[\delta(v-a_2)-K_n(v,a_2)\right].
}
The main part of the proof of \autoref{theo:jac up HE} relies on the following proposition that we shall prove first.
\prop{prop: limit C hat}{For the Jacobi ensemble with P\'olya weight $w_{\rm Jac}(x)=x^\alpha (1-x)^{\beta+n-1} \Theta(1-x)$, $\alpha>-1,\beta>0$ fixed, for all $(r,\lambda_1,\lambda_2)\in \mathbb{R}\times\mathbb{R}_+^2$, the limit
\Eq{}{
\lim_{n\to\infty}\sqrt{n}\,\mathbf{C}_n\left(\frac{1}{2}-\frac{r}{\sqrt{n}};1-\frac{\lambda_1}{n^2},1-\frac{\lambda_2}{n^2} \right)=\frac{4}{\sqrt{\pi }} e^{-4 r^2} \left(\frac{\lambda_1}{\lambda_2} \right)^{\beta/2}J_{\beta }\left(2 \sqrt{\lambda_1}\right)J_{\beta }\left(2 \sqrt{\lambda_2}\right)
}
holds point-wise.}
\begin{proof}[Proof of \autoref{prop: limit C hat}]
Explicitly,
\Eq{eq:cov}{
\mathbf{C}_n(r;a_1,a_2)=&\,\Theta(r-a_2)\int _{0}^{\infty}dt \frac{1}{a_1} \varphi_n\left(\frac{a_1}{r},t\right) K_n\left(a_1,-rt\right)\\
&- \int _{0}^{\infty}dt \int_0^r\frac{dv}{v} \varphi_n\left(\frac{v}{r},t\right) K_n\left(a_1,-rt\right)K_n(v,a_2).
}
From there, one first notices that, after setting $a_2=1-\tfrac{\lambda_2}{n^2}$, the Heaviside function in \eqref{eq:cov} will vanish in the limit as for all fixed $(r,\lambda_2) \in\, \mathbb{R} \times \mathbb{R}_+$,
\Eq{}{
\lim_{n\to \infty}\Theta\left(\frac{r}{\sqrt{n}}+\frac{\lambda_2}{n^2}-\frac{1}{2}\right)=0.
}
The term of interest is therefore
\Eq{eq: scaled TK}{
\int_0^{r}dv \int _{0}^{\infty}dt& K_n\left(1-\tfrac{\lambda_1}{n^2},-rt\right) K_n\left(v,1-\tfrac{\lambda_2}{n^2}\right)\frac{\left(1-\frac{v}{r}\right)^{n-2}}{r^2 (1+t)^{n+2}} \\
&\times\big[(1+rt)-(1-v)-\tfrac{1}{n}(r+vt) \big].
}
 
\rem{rem: dom r}{
We keep $r$ as it is, at the moment, because replacing it by $\frac{1}{2}-\frac{r}{\sqrt{n}}$ does not yield any simplification. We will proceed with the shift and scaling of $r$ at the end, to keep expressions as simple as possible. However, in order not to run in unnecessary difficulty, one needs to take $r$ in compact subsets of $]0,1[$ and can arbitrarily choose to take $r\in \left[\frac{1}{4},\frac{3}{4}\right]$.
}

As we focus on the first term of the asymptotic expansion of \eqref{eq: scaled TK} when $n\to \infty$, it appears that one needs to compute the following integrals
\Eq{eq: block1}{
\int _{0}^{\infty}dt\,   K_n\left(1-\tfrac{\lambda_1}{n^2},-rt\right)(1+t)^{-(n+2)}(1+rt)^\gamma ,\quad \gamma=0,1
}
and
\Eq{eq: block2}{
\int_0^{r}dv\, K_n\left(v,1-\tfrac{\lambda_2}{n^2}\right) \left(1-\frac{v}{r}\right)^{n-2}(1-v)^\gamma, \quad \gamma=0,1.
}
Unfortunately, none of the factors inside the integral admits an individual limit, the goal is therefore to get a good asymptotic of those integrals using Laplace's approximation method. 

%-------------------------------------------------------------------------------------------------------
\subsubsection*{Asymptotic expansion of the kernel}

The first step is to get an asymptotic expansion, as $n\to\infty$, of the integrands in \eqref{eq: block1} and \eqref{eq: block2}, to be able to use Laplace's method. For this, we need the following propositions.
From \cite[p.196]{Szegoe1975} we have the following asymptotic for Jacobi polynomials outside the orthogonality interval.
\prop{prop:asymp Jac1}{
For all $x\in \mathbb{C}\setminus[0,1]$ and arbitrary fixed $\alpha,\beta \in \mathbb{R}$ an asymptotic expansion of $P_{n}^{(\alpha,\beta)}$ is
\Eq{eq:asymp Jac1}{
P_{n}^{(\alpha,\beta)}(1-2x)=A_{n}^{(\alpha,\beta)}(x)\left[1+O\left(\frac{1}{n}\right)\right],\quad \text{as}\quad n\to\infty,
}
where 
\Eq{}{
A_{n}^{(\alpha,\beta)}(x):=\frac{1}{2 \sqrt{\pi  n}}\frac{\left(1+\sqrt{\frac{x}{x-1}}\right)^{\alpha +\beta +n+1}}{ \left(1-\sqrt{\frac{x}{x-1}}\right)^n \left(\frac{x}{x-1}\right)^{\frac{\alpha }{2}+\frac{1}{4}}}.
}
The expansion holds uniformly in the exterior of any closed contour enclosing $[0,1]$, in the sense the ratio goes to $1$.
}
The Mehler–Heine formula \cite{Szegoe1975} gives the following asymptotic for Jacobi polynomials near the right boundary of the orthogonality interval.
\prop{prop:asymp Jac2}{
For any fixed $\alpha>-1,\beta>0$, an asymptotic expansion of $P_{n}^{(\alpha,\beta)}$ is
\Eq{}{
P_{n}^{(\alpha,\beta)}\left(\frac{2\lambda}{n^2}-1 \right) =B_{n}^{(\alpha,\beta)}(\lambda)+O\left(n^{\beta-1}\right),\quad \text{as}\quad n\to\infty,
}
where 
\Eq{}{
B_{n}^{(\alpha,\beta)}(\lambda):=(-1)^n n^\beta \lambda^{-\frac{\beta }{2}} J_{\beta }\left( 2\sqrt{\lambda}\right).
}
The expansion holds uniformly for $\lambda$ in compact subsets of $\mathbb{R}_+$.
}
This leads to the following corollaries, where we retain the dependence of one of the kernel's arguments in the error term. Namely, the argument on which we integrate; cf. \eqref{eq: block1} and \eqref{eq: block2}.
\coro{coro:Asymp kernel1}{
For any fixed , $\alpha>-1,\beta>0$, the expansion
\Eq{eq:Asymp kernel1}{K_n\left(1-\frac{\lambda}{n^2},-x\right)\underset{n\to \infty}{=} -n  \left(\frac{\lambda}{n^2} \right)^\beta \left[ B_{n}^{(\alpha,\beta)}(\lambda)\frac{P_{n}^{(\alpha,\beta-1)}(1+2x)}{1+x}+O\left(n^{\beta-1}\frac{\mathcal{E}_n(-x)}{(1+x)}\right)\right]
}
holds uniformly for $\lambda$ and $x$ in compact subsets of $\mathbb{R}_+$ and $\mathcal{E}_n$ is the polynomial
\Eq{eq: error term jac}{
\mathcal{E}_n(x):=P^{(\alpha,\beta-1)}_{n}\left(1-2x\right)+P^{(\alpha,\beta)}_{n}\left(1-2x\right)+ P^{(\alpha,\beta)}_{n-1}\left(1-2x\right).
}
}
\rem{rem: error term O}{
Note that $B_{n}^{(\alpha,\beta)}(\lambda)\underset{n\to\infty}{=}O(n^\beta)$, therefore, the big $O$ in \eqref{eq:Asymp kernel1} is really an error term.
}
\begin{proof}[Proof of \autoref{coro:Asymp kernel1}]
Using the Christoffel-Darboux formula \eqref{eq:CD formula} and \autoref{prop:asymp Jac2} we have,
\Eq{}{
K_n\left(1-\frac{\lambda}{n^2},-x\right)
\underset{n\to \infty}{=}& a_n \left(\frac{\lambda}{n^2} \right)^\beta \frac{1}{1+x}\bigg[\left(B_{n}^{(\alpha,\beta)}(\lambda)+O\left(n^{\beta-1}\right)\right) P_{n-1}^{(\alpha,\beta)}(1+2x)\\
&-\left(B_{n-1}^{(\alpha,\beta)}(\lambda)+O\left(n^{\beta-1}\right)\right)P_{n}^{(\alpha,\beta)}(1+2x)\bigg].
}
Note that the respective domains of the polynomials do not intersect. Thus, there is no cancellation of the leading terms. 

First, one has to notice that 
\Eq{eq:bn bn-1}{
B_{n-1}^{(\alpha,\beta)}(\lambda) \underset{n\to \infty}{=} -B_{n}^{(\alpha,\beta)}(\lambda)\left[1+O\left(\frac{\beta}{n}\right)\right].
}
Therefore, one can factorise by $B_{n}^{(\alpha,\beta)}(\lambda)$. Then, using the relation 
\Eq{}{
(2n+\alpha+\beta+1)P^{(\alpha,\beta)}_{n}\left(x\right)=(n+\alpha+\beta+1)P^{(\alpha,\beta+1)}_{n}\left(x\right)+(n+\alpha)P^{(\alpha,\beta+1)}_{n-1}\left(x\right),
}
one gets 
\Eq{eq:pn pn-1}{
P^{(\alpha,\beta)}_{n}\left(1+2x\right)+P^{(\alpha,\beta)}_{n-1}\left(1+2x\right) \underset{n\to \infty}{=} 2 P^{(\alpha,\beta-1)}_{n}\left(1+2x\right)+O\left(\frac{\mathcal{E}_n(-x)}{n}\right),
}
which yields
\Eq{}{K_n\left(1-\frac{\lambda}{n^2},-x\right)\underset{n\to \infty}{=} 2 a_n  \left(\frac{\lambda}{n^2} \right)^\beta \left[B_{n}^{(\alpha,\beta)}(\lambda)\frac{P_{n}^{(\alpha,\beta-1)}(1+2x)}{1+x}+O\left(n^{\beta-1}\frac{\mathcal{E}_n(-x)}{(1+x)}\right)\right].
}
Finally, using the following asymptotic equivalent for $a_n$ \eqref{eq: coeff an},
\Eq{eq}{
a_n \underset{n\to \infty}{\sim} -\frac{n}{2},
}
finishes the proof.
\end{proof}
\coro{coro:Asymp kernel2}{
For any fixed $\alpha>-1,\beta>0$, the expansion
\Eq{eq:Asymp kernel2}{K_n\left(v,1-\frac{\lambda}{n^2}\right)\underset{n\to \infty}{=} -n\, v^\alpha (1-v)^{\beta-1}\left[B_{n}^{(\alpha,\beta)}(\lambda)P_{n}^{(\alpha,\beta-1)}(1-2v)+O\left(n^{\beta-1}\mathcal{E}_n(v)\right)\right].
}
holds uniformly for $\lambda$ in compact subsets of $\mathbb{R}_+$ and $v\in [0,1-\epsilon]$, for any $\epsilon\in]0,1[$. The polynomial $\mathcal{E}_n$ is given by \eqref{eq: error term jac}.
}

\begin{proof}[Proof of \autoref{coro:Asymp kernel2}]
In the same spirit as the proof of  \autoref{coro:Asymp kernel1}, we use the Christoffel-Darboux formula \eqref{eq:CD formula} and \autoref{prop:asymp Jac2}. This yields
\Eq{}{
K_n\left(v,1-\frac{\lambda}{n^2}\right)
\underset{n\to \infty}{=}& a_n v^\alpha (1-v)^\beta 
\bigg[P_{n}^{(\alpha,\beta)}(1-2v)\left(B_{n-1}^{(\alpha,\beta)}(\lambda)+O\left(n^{\beta-1}\right)\right) \\
&-P_{n-1}^{(\alpha,\beta)}(1-2v)\left(B_{n}^{(\alpha,\beta)}(\lambda)+O\left(n^{\beta-1}\right)\right)\bigg].
}
Here again, there is no cancellation of the leading terms as the arguments of the kernel are never equal to each other for $n$ large enough, namely $n^2>\lambda/\epsilon$ (cf. \autoref{rem: dom r}). Following the same steps as in the proof of  \autoref{coro:Asymp kernel1}, the result follows.
\end{proof}
From \autoref{coro:Asymp kernel1} and \autoref{coro:Asymp kernel2}, one sees that the problem of getting the large $n$ asymptotic of \eqref{eq: scaled TK} can be broken down to getting the large $n$ asymptotic of the following integrals
\Eq{eq: J int}{
\mathcal{J}^{(\alpha,\beta)}_{n,\gamma}(r):=\int_0^r dv\, v^\alpha (1-v)^{\beta-1+\gamma} P_{n}^{{(\alpha ,\beta-1 )}}(1-2v) \left(1-\frac{v}{r}\right)^{n-2}, \quad \gamma=0,1
}
and
\Eq{eq: I int}{
\mathcal{I}^{(\alpha,\beta)}_{n,\gamma}(r):=\int_0^\infty dt (1+t)^{-(n+2)} \frac{P_{n}^{(\alpha,\beta-1)}(1+2rt)}{(1+rt)^{1-\gamma}},\quad \gamma=0,1.
}
Note that the error terms in \eqref{eq:Asymp kernel1} and \eqref{eq:Asymp kernel2} are composed of the same polynomials as the leading terms, with same degree. Indeed, $O(\mathcal{E}_n(x))$ is some linear combination of the Jacobi polynomials composing $\mathcal{E}_n$. Therefore, after integration, the error terms will respectively be expressed as some linear combination of  $\mathcal{J}^{(\alpha,\beta)}_{n,\gamma}$ and $\mathcal{I}^{(\alpha,\beta)}_{n,\gamma}$, up to a shift in the parameters $\alpha,\beta$, which are $n$-independent. The error terms will remain $1/n$ smaller than the leading terms and can, therefore, be neglected in the limit $n\to\infty$; see \autoref{rem: error term O}.

%-------------------------------------------------------------------------------------------------------
\subsubsection*{Computation of $ \mathcal{J}^{(\alpha,\beta)}_{n,\gamma}(r)$}

Let us focus first on $\mathcal{J}^{(\alpha,\beta)}_{n,\gamma}(r)$, $ \gamma=0,1$, as one does not need to use the Laplace's method in this case. Indeed, one can use Rodrigues' formula for the Jacobi polynomials $P_{n}^{{(\alpha ,\beta)}}$, which, under the mapping $x\mapsto1-2x$, is given by
\Eq{}{
x^\alpha (1-x)^\beta P_{n}^{{(\alpha ,\beta)}}(1-2x)=\frac{1}{n!}\frac{d^n}{dx^n}\left[x^{\alpha+n} (1-x)^{\beta+n}\right].
}
Then, integrating by part, respectively $n-2$ and $n-1$ times, yields
\Eq{}{
\mathcal{J}^{(\alpha,\beta)}_{n,0}(r)=\frac{r^{\alpha +1} (1-r)^{n+\beta-2} }{(n-1)}\left[1-2r+\frac{\alpha-r(\alpha+\beta-1)}{n}\right]
}
and
\Eq{}{
\mathcal{J}^{(\alpha,\beta)}_{n,1}(r)=\frac{r^{\alpha +1} (1-r)^{n+\beta-1} }{(n-1)}\left[1-r+\frac{\alpha-r(\alpha+\beta)}{n}\right].
}

%-------------------------------------------------------------------------------------------------------
\subsubsection*{Asymptotic equivalent of $\mathcal{I}^{(\alpha,\beta)}_{n,\gamma}(r)$}
Let us now find asymptotic equivalents of $\mathcal{I}^{(\alpha,\beta)}_{n,\gamma}(r)$, $\gamma=0,1$. Once again the integral can be computed explicitly, however this yield a sum of hypergeometric functions for which an asymptotic equivalent is hard to get. Instead, we will use Laplace's approximation method. For this, we need the asymptotic for $P_{n}^{(\alpha,\beta-1)}$ given by \autoref{prop:asymp Jac1}. However, one needs to stay away from $t=0$, as the approximation is not uniform there. We then introduce a $n$-independent cut-off $t_\varepsilon>0$ that we can conveniently express in terms of a small enough $\varepsilon>0$ as
\Eq{}{
\varepsilon=\sqrt{\frac{rt_\varepsilon}{1+rt_\varepsilon}}.
}
We will see how small $\varepsilon$ needs to be later on. To estimate the contribution of the integral
\Eq{eq: cut integ}{
\int_0^{t_\varepsilon} dt (1+t)^{-(n+2)} \frac{P_{n}^{(\alpha,\beta-1)}(1+2rt)}{(1+rt)^{1-\gamma}},\quad \gamma=0,1,
}
one can bound $P_{n}^{(\alpha,\beta-1)}(1+2rt)$ by its maximum, which is attained at $t=t_\varepsilon$ and take the very rough upper bound
\Eq{}{
\int_0^{t_\varepsilon} dt (1+t)^{-(n+2)} \frac{P_{n}^{(\alpha,\beta-1)}(1+2rt)}{(1+rt)^{1-\gamma}}\leq \frac{1}{n+1}P_{n}^{(\alpha,\beta-1)}(1+2rt_\varepsilon),\quad \gamma=0,1.
}
To estimate the order of the upper bound as $n\to\infty$, one can use the approximation \eqref{eq:asymp Jac1} as $t_\varepsilon$ is bounded away from $0$ and $n$-independent. This yields, in terms of $\epsilon$,
\Eq{eq: order remainder}{
\frac{1}{n+1}P_{n}^{(\alpha,\beta-1)}(1+2rt_\varepsilon)\underset{n\to\infty}{=}O\left(\frac{1}{n^{3/2}}\left(\frac{1+\varepsilon}{1-\varepsilon}\right)^n\right),
}
for $r$ in compact subsets of $]0,1[$; see \autoref{rem: dom r}. Using this estimate, we will see later that the integral \eqref{eq: cut integ} is negligible compared to its complement.

The plan is thus to apply Laplace's method on the following integral for $n\to \infty$, 
\Eq{}{
\mathcal{Y}^{(\alpha,\beta)}_{n,\gamma}(r):=\int_{t_\varepsilon}^\infty dt (1+t)^{-(n+2)}  \frac{A_{n}^{(\alpha,\beta-1)}(-rt)}{(1+rt)^{1-\gamma}},\quad \gamma=0,1.
}
As the asymptotic expansion \eqref{eq:asymp Jac1} is uniform on the integration domain, we have 
\Eq{eq: error order I}{
\mathcal{I}^{(\alpha,\beta)}_{n,\gamma}(r) \underset{n\to \infty}{=} \mathcal{Y}^{(\alpha,\beta)}_{n,\gamma}(r)\left[1+O\left(\frac{1}{n}\right)\right]+O\left(\frac{1}{n^{3/2}}\left(\frac{1+\varepsilon}{1-\varepsilon}\right)^n\right),
}
for all $r$ in compact subsets of $]0,1[$, where both $\mathcal{I}^{(\alpha,\beta)}_{n,\gamma}(r)$ and $\mathcal{Y}^{(\alpha,\beta)}_{n,\gamma}(r)$, are convergent integrals due to the factor $(1+t)^{-(n+2)}$. Proceeding with a change of variable $y=\sqrt{\frac{rt}{1+rt}}$, the change of measure is then $dt=\frac{2y dy}{r(1-y^2)^2}$ and the integral reads
\Eq{}{
\mathcal{Y}^{(\alpha,\beta)}_{n,\gamma}(r)= \frac{2}{r}\int_{\varepsilon}^1 dy\, y(1-y^2)^{n+1-\gamma} \left(1+y^2(\tfrac{1}{r}-1)\right)^{-(n+2)}  A_{n}^{(\alpha,\beta-1)}\left(-\frac{y^2}{1-y^2}\right).
}
One can then rewrite the integrand
\Eq{}{
\frac{g_{r,\gamma}(y)}{\sqrt{n}}\exp(-n f_r(y))=&\,\frac{2}{r}  y(1-y^2)^{n+1-\gamma} \left(1+y^2(\tfrac{1}{r}-1)\right)^{-(n+2)}  A_{n}^{(\alpha,\beta-1)}\left(-\frac{y^2}{1-y^2}\right)\\
=&\, \frac{1}{r\sqrt{n\pi}}y^{\frac{1}{2}-\alpha}(1-y)^{1-\gamma}\frac{(1+y)^{2n+\alpha+\beta+1-\gamma}}{\left(1+y^2(\tfrac{1}{r}-1)\right)^{n+2} },
}
defining, as follows, the functions $f_r$ and $g_{r,\gamma}$, for which we already give the explicit expressions of the needed derivatives
\Eq{}{
&f_r(y):=\ln(1+y^2(\tfrac{1}{r}-1))-2\ln(1+y),\\
&\\
&\displaystyle f_r'(y)=2\frac{y(1-r)-r}{(1+y)(r+y^2(1-r))},\\
&\\
&\displaystyle f_r''(y)=2\frac{r(1+y)(1+y+4y^2)-2r^2y(1+y)^2-y^2(1+2y)}{(1+y)^2(r+y^2(1-r))^2}\\
}
and 
\Eq{}{
g_{r,\gamma}(y):=&\displaystyle\frac{1}{r\sqrt{\pi}}y^{\frac{1}{2}-\alpha}(1-y)^{1-\gamma}\frac{(1+y)^{\alpha+\beta+1-\gamma}}{\left(1+y^2(\tfrac{1}{r}-1)\right)^{2} },\\
&\\
 g_{r,\gamma}'(y)=&\displaystyle \frac{y^{-\alpha -\frac{1}{2}} (1+y)^{\alpha +\beta }\left(1-y^2\right)^{-\gamma }}{2 r^2 \sqrt{ \pi }   \left(1+\left(\frac{1}{r}-1\right) y^2\right)^3} \bigg[y^2 \left(-2 \alpha +y^2 (-2 \beta +4 \gamma +3)+2 y (\alpha +\beta )-7\right)\\
 &+r \left(y^2-1\right) \left(2 \alpha +y^2 (2 \beta -4 \gamma -3)-2 y (\alpha +\beta )-1\right)\bigg].\\
}
The function $f_r$ admits a unique critical point; corresponding to a minimum located at $y=\frac{r}{1-r}$. For $r>\frac{1}{2}$, the critical point escapes the domain of integration, the minimum is then on the boundary at $y=1$. However, $g_{r,0}(1)=0$, hence the need to expand $g_{r,\gamma}$ to first order and $f_r$ to second order. Note that, despite the term $\left(1-y^2\right)^{-1 }$, $g_{r,1}'$ does not have a singularity at $y=1$, as it is removable.
Denoting, 
\Eq{}{
\widehat{\mathcal{Y}}^{(\alpha,\beta)}_{n,\gamma}(r) := &
\int_\varepsilon^1 \frac{dy}{\sqrt{n}} \left[ g_{r,\gamma}(y_r)+(y-y_r) g_{r,\gamma}'(y_r)\right]\\
\times & \exp(-n\left[f_r(y_r)+(y-y_r) f_r'(y_r)+\frac{(y-y_r)^2}{2} f_r''(y_r) \right]), \quad \gamma=0,1,
}
with 
\Eq{}{
y_r= 
\begin{cases}
    &\displaystyle \frac{r}{1-r}, \quad 0< r \leq \frac{1}{2}\\
    &\\
    &\displaystyle 1, \quad \frac{1}{2}< r < 1\\
\end{cases}.
}
We thus have
\Eq{}{
 \mathcal{Y}^{(\alpha,\beta)}_{n,\gamma}(r) \underset{n\to \infty}{=} 
\widehat{\mathcal{Y}}^{(\alpha,\beta)}_{n,\gamma}(r)\left[1+O\left(\frac{1}{\sqrt{n}}\right)\right],
}
uniformly for $r$ in compact subsets of $]0,1[$. Note that the relative error is, here, of order $n^{-1/2}$ and not $n^{-1}$ due to $g_{r,0}$ vanishing at $y=1$; see \cite[p.37]{Temme2014}. To compute explicitly $\widehat{\mathcal{Y}}^{(\alpha,\beta)}_{n,\gamma}(r)$, one proceeds with the change of variable $u=y_r-y$. The expression becomes
\Eq{}{ 
\widehat{\mathcal{Y}}^{(\alpha,\beta)}_{n,\gamma}(r)=\int_{y_r-1}^{y_r-\varepsilon} \frac{du}{\sqrt{n}} &\left[ g_{r,\gamma}(y_r)-u g_{r,\gamma}'(y_r)\right] \exp(-n\left[f_r(y_r)-u f_r'(y_r)+\frac{u^2}{2} f_r''(y_r) \right])
}
and can be computed using the following lemma, which is obtained with a simple integration by part.
\lem{lem:erf integ}{ Let $a,b,c,d,h,x_0,x_1 \in \mathbb{R}$, $a\neq 0$,
\Eq{}{
&\int_{x_0}^{x_1} (d-h x) \exp \left(-\tfrac{a}{2}x^2 -b x+c \right) \, dx\\
=&\,\tfrac{h}{a} \left(e^{-\tfrac{a x_1^2}{2}-b x_1+c}-e^{-\tfrac{a x_0^2}{2}-b x_0+c}\right)+\sqrt{\tfrac{\pi }{2}} e^{\tfrac{b^2}{2 a}+c}\tfrac{ (a d+b h) }{a^{3/2}} \left(\erf\left(\tfrac{a x_1+b}{\sqrt{2 a}}\right)-\erf\left(\tfrac{a x_0+b}{\sqrt{2 a}}\right)\right),
}
with 
\Eq{}{
\erf(x)=\frac {2}{\sqrt {\pi }}\int _{0}^{x}e^{-t^{2}}\,\mathrm {d} t, \quad x\in \mathbb{R}.
}
}
Using \autoref{lem:erf integ} and the following expressions
\Eq{}{
f_r'(1)=(1-2r),\quad f_r''(1)=-\left(4r^2-6r+\frac{3}{2}\right),\quad f_r''\left(\frac{r}{1-r}\right)=2\frac{(1-r)^3}{r} 
}
and
\Eq{}{
 g_{r,\gamma}'\left(\frac{r}{1-r}\right)=-\frac{4 r^2 (\alpha +\beta -\gamma -3)-r (6 \alpha +2 \beta -10)+2 \alpha-1}{2\sqrt{\pi  }\, r^{\alpha +\frac{3}{2}} \left(1-r\right)^{\beta -\frac{1}{2}}}\left(\frac{(1-r)^2}{1-2 r}\right)^{\gamma },
}
\Eq{}{
g_{r,\gamma}'( 1)=-\frac{ 2^{\alpha +\beta -1}}{\sqrt{\pi  }} r\bigg(4+\gamma(\alpha-\beta+3-8r)\bigg),
}
one gets, after simplifications, for $r\leq \frac{1}{2}$,
\Eq{eq: approx below}{
\widehat{\mathcal{Y}}^{(\alpha,\beta)}_{n,\gamma}(r)  \underset{n\to \infty}{=}& \, \frac{(1-2 r)^{-\gamma }(1-r)^{2\gamma}}{4n \sqrt{\pi }\, r^{\alpha+1} (1-r)^{\beta +n+2}}\Bigg[ 2 \sqrt{ \pi }  (1-2 r) r\\
&\times \bigg[\erf\left(\sqrt{n} (1-2 r) \sqrt{\tfrac{1-r}{r}}\right)+\erf\left(\sqrt{n} (r-(1-r) \epsilon)\sqrt{\tfrac{1-r}{r}}\right)\bigg]\\
&+\frac{1}{\sqrt{n}}\sqrt{\tfrac{r}{1-r}}\left(e^{-n (1-2 r)^2 \frac{(1-r)}{r}}-e^{-n(r-(1-r) \epsilon )^2\frac{ (1-r) }{r}}\right)\\
&\times\left(4 r^2 (\alpha +\beta -\gamma-3)-2 r (3 \alpha +\beta -5)+2 \alpha -1\right)\Bigg]\left(1+O\left(\frac{1}{\sqrt{n}}\right)\right)
}
and, for $r> \frac{1}{2}$, 
\Eq{}{
\widehat{\mathcal{Y}}^{(\alpha,\beta)}_{n,\gamma}(r)  \underset{n\to \infty}{=} &\, \frac{2^{2n+\alpha +\beta+2(1-\gamma)} r^{n+1}(\alpha-\beta+7-8r)^\gamma}{n\sqrt{\pi }  \left(4 (3-2 r) r-3\right)^{\frac{3}{2}}} \Bigg[\frac{1}{\sqrt{n}}\sqrt{4 (3-2 r) r-3} \\
&+\sqrt{\pi}\left[ 2 r-1-\gamma \left(\frac{4 (3-2 r) r-3}{\alpha-\beta+7-8r}\right) \right]\exp(\frac{n (2 r-1)^2}{4 (3-2 r) r-3})\\
&\times  \left(\erf\left(\frac{\sqrt{n} (2 r-1)}{\sqrt{4 (3-2 r) r-3}}\right)-1\right)\Bigg]\left(1+O\left(\frac{1}{\sqrt{n}}\right)\right).
}
Note that the $\varepsilon$ dependence completely drops out here. This is also the case in \eqref{eq: approx below} if one chooses $\varepsilon<\frac{r}{1-r}$. As $r$ is bounded away from $0$, this is certainly always possible; see \autoref{rem: dom r}. Actually, coming back to the remainder integral \eqref{eq: cut integ}, one needs to take $\varepsilon$ slightly smaller, to keep its relative contribution negligible. Indeed, taking $\varepsilon<\frac{r}{1-r}$, for $r$ in compact subsets of $]0,1[$,
\Eq{}{
\widehat{\mathcal{Y}}^{(\alpha,\beta)}_{n,\gamma}(r)  \underset{n\to \infty}{=}O\left(\frac{1}{n(1-r)^n}\right),\quad r\leq 1/2,
}
and
\Eq{}{
\widehat{\mathcal{Y}}^{(\alpha,\beta)}_{n,\gamma}(r)  \underset{n\to \infty}{=}O\left(\frac{(4r)^n}{n}\right), \quad r>1/2.
}
The latter is obtained using the asymptotic expansion of the complementary error function
\Eq{eq: asymp erfc}{
1-\erf(x)&\underset{x\to \infty}{=} \frac {e^{-x^{2}}}{x{\sqrt {\pi }}}\left[1-\frac{1}{2 x^2}+O\left(\frac{1}{x^4} \right) \right].
}
Therefore, comparing the order of $\widehat{\mathcal{Y}}^{(\alpha,\beta)}_{n,\gamma}(r)$ with the order of the remainder integral \eqref{eq: order remainder}, one needs to choose $\varepsilon$ such that, for $r$ on compact subsets of $]0,1[$,
\Eq{}{
\begin{cases}
    &\left(\frac{1+\varepsilon}{1-\varepsilon}\right)(1-r)<1,\quad r\leq1/2\\
    &\left(\frac{1+\varepsilon}{1-\varepsilon}\right)\frac{1}{4r}<1,\quad r>1/2
\end{cases} ,
}
which can be reduced to 
\Eq{}{
\varepsilon<\min\left\{ \frac{r}{2-r},\frac{1}{3}\right\}.
}
This choice of $\varepsilon$ guarantees that the remainder integral \eqref{eq: cut integ} is exponentially small compared to $\mathcal{Y}^{(\alpha,\beta)}_{n,\gamma}(r)$ as $n\to\infty$. Consequently, the expansion \eqref{eq: error order I} reduces to
\Eq{eq: error order I 2}{
\mathcal{I}^{(\alpha,\beta)}_{n,\gamma}(r) \underset{n\to \infty}{=} \mathcal{Y}^{(\alpha,\beta)}_{n,\gamma}(r)\left[1+O\left(\frac{1}{n}\right)\right]\underset{n\to \infty}{=} 
\widehat{\mathcal{Y}}^{(\alpha,\beta)}_{n,\gamma}(r)\left[1+O\left(\frac{1}{\sqrt{n}}\right)\right].
}

%----------------------------------------------------------------------------------------
\subsubsection*{Asymptotic equivalent of the $1,k$-cross-covariance function}

With all the asymptotic expansions found above, we are in good position to find the desired asymptotic equivalent of the $1,k$-cross-covariance density. Let us first notice that
\Eq{}{
\int_0^r dv\, v^{\alpha+1} (1-v)^{\beta-1} P_{n}^{{(\alpha ,\beta-1 )}}(1-2v) \left(1-\frac{v}{r}\right)^{n-2}=\mathcal{J}^{(\alpha,\beta)}_{n,0}(r)-\mathcal{J}^{(\alpha,\beta)}_{n,1}(r)
}
and
\Eq{}{
\int_0^\infty dt (1+t)^{-(n+2)} \frac{P_{n}^{(\alpha,\beta-1)}(1+2rt)}{1+rt} t=\frac{\mathcal{I}^{(\alpha,\beta)}_{n,1}(r)-\mathcal{I}^{(\alpha,\beta)}_{n,0}(r)}{r}.
}
Then, gathering the asymptotic expansions \eqref{eq:Asymp kernel1} and \eqref{eq:Asymp kernel2}, along with the integrals $\mathcal{J}^{(\alpha,\beta)}_{n,\gamma}(r)$ \eqref{eq: J int} and $\mathcal{I}^{(\alpha,\beta)}_{n,\gamma}(r)$ \eqref{eq: I int}, one gets the following asymptotic equivalent
\Eq{eq: asymp interm}{
&\int_0^{r}dv \int _{0}^{\infty}dt K_n\left(1-\tfrac{\lambda_1}{n^2},-rt\right) K_n(v,1-\tfrac{\lambda_2}{n^2}) \left(1-\frac{v}{r}\right)^{n-2} (1+t)^{-(n+2)}\\
&\times \frac{1}{r^2}\big[(1+rt)-(1-v)-\tfrac{1}{n}(r+vt) \big]\\
\underset{n\to \infty}{\sim}& n^{2-2\beta} \lambda_1^{\beta} B_{n}^{(\alpha,\beta)}(\lambda_1)B_{n}^{(\alpha,\beta)}(\lambda_2) \mathcal{G}^{(\alpha,\beta)}_{n}(r)
}
with
\Eq{}{
\mathcal{G}^{(\alpha,\beta)}_{n}(r):=& \frac{1}{r^2} \bigg( \mathcal{I}^{(\alpha,\beta)}_{n,0}(r)\left[\left(\frac{1}{nr}-\frac{r}{n} \right)\mathcal{J}^{(\alpha,\beta)}_{n,0}(r)-\left(1+\frac{1}{nr} \right)\mathcal{J}^{(\alpha,\beta)}_{n,1}(r)\right] \\
&+ \mathcal{I}^{(\alpha,\beta)}_{n,1}(r)\left[\left(1-\frac{1}{nr} \right)\mathcal{J}^{(\alpha,\beta)}_{n,0}(r)+\frac{1}{nr}\mathcal{J}^{(\alpha,\beta)}_{n,1}(r)\right]  \bigg).\\
}
The explicit expression of $\mathcal{G}^{(\alpha,\beta)}_{n}(r)$ is rather cumbersome, however, there is some simplification by noting
\Eq{}{
&\left(\frac{1}{nr}-\frac{r}{n} \right)\mathcal{J}^{(\alpha,\beta)}_{n,0}(r)-\left(1+\frac{1}{nr} \right)\mathcal{J}^{(\alpha,\beta)}_{n,1}(r)\\
=&-\frac{r^{\alpha +1} (1-r)^{\beta +n-1}}{n} \left(1-r+\frac{1+\alpha -r (\alpha +\beta -1)}{n}\right)
}
and
\Eq{}{
&\left(1-\frac{1}{nr} \right)\mathcal{J}^{(\alpha,\beta)}_{n,0}(r)+\frac{1}{nr}\mathcal{J}^{(\alpha,\beta)}_{n,1}(r)\\
=&\frac{r^{\alpha +1} (1-r)^{\beta +n-2}}{n} \left(1-2r+\frac{1+\alpha -r (\alpha +\beta )}{n}\right).
}
The asymptotic equivalent $\widehat{\mathcal{Y}}^{(\alpha,\beta)}_{n,\gamma}(r)$ of the integrals $\mathcal{J}^{(\alpha,\beta)}_{n,\gamma}(r)$ are, at the moment, defined piecewise, depending on the value of $r$. Both pieces must agree at $r=1/2$, which is precisely the value around which we are zooming.

For $r\leq 1/2$, using the asymptotic equivalent $\widehat{\mathcal{Y}}^{(\alpha,\beta)}_{n,\gamma}(r)$ and proceeding with simplifications, one gets the following asymptotic equivalent 
\Eq{}{
\mathcal{G}^{(\alpha,\beta)}_{n}(r) \underset{n\to \infty}{\sim}&\frac{-1}{ n^2 (1-r)^3 \sqrt{\pi }}\Bigg[ \frac{1}{\sqrt{n}} \left( \sqrt{r(1-r)}-\frac{1}{n}\frac{P_0(r)}{4 \sqrt{r(1-r)}}\right)e^{-n\frac{ (1-r)}{r} (1-2 r)^2}\\
&+\frac{\sqrt{\pi }}{2 n} \left[1+\erf\left((1-2 r) \sqrt{\frac{n (1-r)}{r}}\right)\right] (r (\alpha +\beta -2)-\alpha )\Bigg]
}
with the polynomial
\Eq{}{
P_0(r):=&\,2 r^2 \left[\alpha ^2+2 \alpha  (\beta -2)+(\beta -4) \beta +6\right]\\
&+r \left[\alpha  (-4 \alpha -4 \beta +9)+\beta -6\right]+\alpha  (2 \alpha -1).
}
From there, one can proceed with the change of variable $r \mapsto\frac{1}{2}-\frac{r}{\sqrt{n}}$. This gives 
\Eq{}{
\mathcal{G}^{(\alpha,\beta)}_{n}\left(\frac{1}{2}-\frac{r}{\sqrt{n}}\right) \underset{n\to \infty}{=} -\frac{4 e^{-4 r^2}}{\sqrt{\pi } n^{5/2}} +O\left(\frac{1}{n^3} \right).
}
For $r> 1/2$, one gets the following asymptotic equivalent 
\Eq{}{
\mathcal{G}^{(\alpha,\beta)}_{n}(r) \underset{n\to \infty}{\sim}&-\frac{2^{\alpha +\beta +2 n} r^{\alpha +n} (1-r)^{\beta +n-2}}{n^{5/2} \left(8 r^2-12 r+3\right)^2 \sqrt{\pi }} \Bigg[\left(8 r^2-12 r+3\right)\left( P_1(r)+\frac{1}{n}P_2(r)\right)\\
&-\sqrt{\pi } \sqrt{-n \left(8 r^2-12 r+3\right)} \left( (1-2 r)^2 (-\alpha +\beta +2 r)+\frac{1}{n}P_3(r) \right)\\
&\times e^{-\frac{n (1-2 r)^2}{8 r^2-12 r+3}}\left(\erf\left(\frac{n (2 r-1)}{\sqrt{-n \left(8 r^2-12 r+3\right)}}\right)-1\right)\Bigg],
}
with
\Eq{}{
P_1(r):=&\,12 r^2-2 r (\alpha -\beta +7)+\alpha -\beta +3,\\
P_2(r):=&\,4 r^2 (\alpha +\beta +1)-r \left(\alpha ^2+7 \alpha -\beta ^2+3 \beta +8\right)+(\alpha +1) (\alpha -\beta +3),\\
P_3(r):=&\, 8 r^3+2 r^2 \left(-\alpha ^2+\alpha +\beta ^2+\beta -6\right)\\
&+r ((\alpha -\beta ) (3 \alpha +\beta )-2 \beta +2)-(\alpha +1) (\alpha -\beta ).
}
Then, using the asymptotic expansion of the complementary error function \eqref{eq: asymp erfc} and the same change of variable $r \mapsto\frac{1}{2}-\frac{r}{\sqrt{n}}$, one arrives at the following asymptotic equivalent 
\Eq{}{
\mathcal{G}^{(\alpha,\beta)}_{n}\left(\frac{1}{2}-\frac{r}{\sqrt{n}}\right) \underset{n\to \infty}{\sim}& -\frac{4}{\sqrt{\pi } n^{5/2}}\left(1+\frac{2 r}{\sqrt{n}}\right)^n \left(1-\frac{2 r}{\sqrt{n}}\right)^n\bigg[1+\frac{1}{n}P_2(1/2)  \\
&-\frac{1}{\sqrt{n}} \sqrt{\pi}e^{4r^2}\left(\erf(2r)-1\right) \bigg].
}
The last step is to use the expansion
\Eq{}{
 \left(1-\frac{4 r^2}{n}\right)^n\underset{n\to \infty}{=}e^{-4r^2}\left(1+O\left(\frac{1}{n}\right)\right)
}
to arrive once again at
\Eq{}{
\mathcal{G}^{(\alpha,\beta)}_{n}\left(\frac{1}{2}-\frac{r}{\sqrt{n}}\right) \underset{n\to \infty}{=} -\frac{4 e^{-4 r^2}}{\sqrt{\pi } n^{5/2}} +O\left(\frac{1}{n^3} \right).
}
Plugging this back in \eqref{eq: asymp interm} and then in \eqref{eq:cov} finishes the proof.

\end{proof}

\begin{proof}[Proof of \autoref{theo:jac up HE}]
Starting from \autoref{theo:poly ensemble} and following the same arguments as the proof of \autoref{coro: cov poly ensemble}, one gets 
\Eq{}{
\mathrm{cov}_{1,k}(r;a_1,\ldots,a_k)=&\frac{(n-k)!}{(n-1)!}\sum_{l,m=1}^k(-1)^{l+m}\det[K_n(a_b,a_c)]_{\substack{b\in\llbracket 1,k\rrbracket\setminus\{l\}\\c\in\llbracket 1,k\rrbracket\setminus\{m\}}}\mathbf{C}_n(r;a_l,a_m).
}
Proceeding with the change of variables $r\to \frac{1}{2}-\frac{r}{\sqrt{n}}$ and $a_j\to 1-\frac{a_j}{n^2}$, $j=1,\ldots,k$, and scaling by a factor $n^{k+2}$, one gets
\Eq{}{
&\frac{n^{k+2}}{n^{2k+1/2}}\mathrm{cov}_{1,k}\left(\frac{1}{2}-\frac{r}{\sqrt{n}}; 1-\frac{a_1}{n^2},\ldots,1-\frac{a_k}{n^2} \right)\\
=&\,\frac{(n-k)!\,n^k}{n!}\sum_{l,m=1}^k \sqrt{n} \mathbf{C}_n\left(\frac{1}{2}-\frac{r}{\sqrt{n}};1-\frac{a_l}{n^2},1-\frac{a_m}{n^2}\right) \\
&\qquad\qquad\qquad\times(-1)^{l+m}\det[\frac{1}{n^2}K_n\left(1-\frac{a_b}{n^2},1-\frac{a_c}{n^2}\right)]_{\substack{b\in\llbracket 1,k\rrbracket\setminus\{l\}\\c\in\llbracket 1,k\rrbracket\setminus\{m\}}}.
}
Then, taking the limit $n\to \infty$ with the help of \eqref{eq: lim factorial}, the limiting kernel \eqref{eq: limiting kernel Jac UHE} and \autoref{prop: limit C hat}, one gets, using Jacobi's formula (cf. \cite[Eq.(4.10)]{Allard2025b}), that $\mathrm{cov}_{1,k\, {\rm SH}}^{(\infty)}(r; a_1,\ldots,a_k)$ is equal to
\Eq{eq: diff det form}{
  \frac{4}{\sqrt{\pi }} e^{-4 r^2} \partial_\mu \det\left[ \left(\frac{a_b}{a_c} \right)^{\beta/2}4 K_{\beta}^{\mathrm{Bes}}(4a_b,4a_c)+\mu \left(\frac{a_b}{a_c} \right)^{\beta/2} J_{\beta}(2\sqrt{a_b }) J_{\beta}(2\sqrt{a_c }) \right]_{b,c=1}^k \Bigg|_{\mu=0}.\\
}
The factor $\left(\frac{a_b}{a_c} \right)^{\beta/2}$ can be eliminated by elementary transformations on the determinant.
Alternatively, as the dependence of $a_b$, $a_c$ factorizes in the scaling limit of $\mathbf{C}_n$, \eqref{eq: diff det form} can also be put in the form 
\Eq{}{
-\frac{4^k}{\sqrt{\pi }} e^{-4 r^2}\det\left(\begin{array}{c c} 
    	0\quad &  J_{\beta}(2\sqrt{a_c })  \\
     J_{\beta}(2\sqrt{a_b }) \quad & K_{\beta}^{\mathrm{Bes}}(4a_b,4a_c) 
\end{array}\right)_{b,c=1}^k.\\
}

\end{proof}

%=================================================================
\section{Discussion}\label{Discussion}

In the first part of this article, we have exploited the new results of \cite{Allard2025b}, for polynomial ensembles, to study their extension to the large $n$ limit, which are thus also new. First, we focused on the double scaling limit around the origin for the $1,k$-point correlation function between one eigenradius and $k$ singular values as well as the $1,k$-cross-covariance density function, see \autoref{theo: 1,k poly ensemble}, \autoref{coro: cov poly ensemble} and \autoref{coro:polya ensemble}. This was motivated by the fact all ensembles potentially have a hard edge at the origin for both singular values and eigenradii, as they are non-negative. An interesting finding is the existence of a universal (at least for polynomial ensembles) scaling ratio of $1/n$ between the scale of the smallest squared singular value and scale of the smallest squared eigenradius.

Under  \autoref{assum:exist HE PE} guaranteeing the existence of the limits, the formula for the $1,k$-point correlation function involves the double scaling limit of the kernel of the corresponding determinantal point process on the singular values. The result makes clear that all ensembles for which the limiting kernel are the same will share the same limiting $1,k$-point correlation function around the origin and same limiting $1,k$-cross-covariance function (cf. \eqref{eq:1,kpt poly inf}).  

For Pólya ensembles, this translates into having same limiting Pólya weight function. As for the $1,k$-point function, it admits a more explicit form given by \autoref{coro:polya ensemble}. We give the example of three classical Pólya ensembles; Laguerre, Jacobi and Cauchy-Lorentz ensembles, sharing the same limiting weight (cf. \autoref{ex:classic polya}). To further illustrate this, \autoref{fig: limiting cov} shows the similarity between the plots of the, appropriately rescaled, cross-covariance for Jacobi and Laguerre ensembles at $n=25$ and the plot of their shared limiting cross-covariance. The plots of \autoref{fig: limiting cov} also show how the deterministic constraints from Weyl's inequalities~\cite{Weyl1949} translate on a probabilistic level and survive the large $n$ limit. Here, the plot of the cross-covariance function around the origin is reminiscent of the fact that the smallest eigenradius is bounded from below by the smallest singular value. 

We then give the example of Muttalib-Borodin ensembles to emphasize how the main results \autoref{coro:polya ensemble} and \autoref{theo: 1,k polya ensemble} can be applied. 

In a second part, we looked at the particular case of Jacobi ensembles which enjoy a second hard edge for the singular values. It turns out that the corresponding edge for the eigenradii is soft. The formula of the limiting $1,k$-cross-covariance density at the soft-hard edge involves some Bessel function for the singular value part, as expected, and is Gaussian in the eigenradius (cf. \autoref{theo:jac up HE}). We note that this limiting $1,k$-cross-covariance at the soft-hard edge is smaller by a factor $1/n$ than the limiting $1,k$-point function and can therefore be interpreted as a correction term. Indeed, as the eigenradius does not share the same edge as the singular values in this case, the Weyl's inequalities are much less felt and this translates in a vanishing local cross-covariance, in the large $n$ limit. From this perspective, the Gaussian dependence on the eigenradius is to put in relation to a recurring phenomenon for $\beta$-ensembles, as pointed out in \cite{forrester2021}, that the leading correction term (cf.\autoref{rem: speed decrease 2}) at the edge (soft and hard) is related to the limiting distribution--- in this case the complementary error function---by a derivative operation.

Coming back to the main results, we strongly believe the formula \eqref{eq:1,kpt poly inf} for the limiting $1,k$-point function (equivalently, the $1,k$-cross-covariance) at the origin remains true even in the case the limiting kernel is not a function anymore but a general distribution e.g. a weighted Dirac delta function, as the formula  \eqref{eq:1,kpt poly inf} still makes sense in the distributional sense, in this case. However, proving it will require a different proof, as one would not have Lebesgue's dominated convergence theorem anymore, to interchange the limit $n\to \infty$ and the integrals. We also believe that \autoref{assum:exist HE} strictly imply \autoref{assum:exist HE PE}. In \autoref{Appendix}, we show that the technical restriction \eqref{eq: asymp pj polynom} on polynomial ensembles is always verified for P\'olya ensembles verifying \autoref{assum:exist HE}, as well as their compositions; see \autoref{rem: closure}. Some assumptions can possibly be lifted with some further analysis and making use of some property of the underlying Pólya frequency functions such as their log-concavity \cite{Schoenberg1951,Foerster2020,Groechenig2020} as well as the log-convexity of the Mellin transform of probability measures (on each subinterval of $\mathbb{R}_+$ where it is finite). 

While \autoref{theo: 1,k poly ensemble} targets the largest class of ensembles where the formula \eqref{eq:1,kpt poly inf} holds, it relies on \autoref{assum:exist HE PE} which are hard to check in practice as they require finding the bi-orthonormal system of functions composing the kernel of the corresponding determinantal point process and this for all $n\in \mathbb{N}$. On the other hand, \autoref{theo: 1,k polya ensemble} and \autoref{coro:polya ensemble} are much more useful, as \autoref{assum:exist HE} are much easier to check. Indeed, the assumptions bear only on the Pólya weight function, which is easier to find in practice.

This barely diminishes our results as the class of Pólya ensembles, albeit being smaller than the class of polynomial ensembles, is very large and comprises most of the classical ensembles in the field of random matrices along with their compositions due to their closure property, as discussed in \autoref{Introduction} and \autoref{rem: closure}.

\begin{figure}[H]
        \centering
        \includegraphics[width=1.1\textwidth]{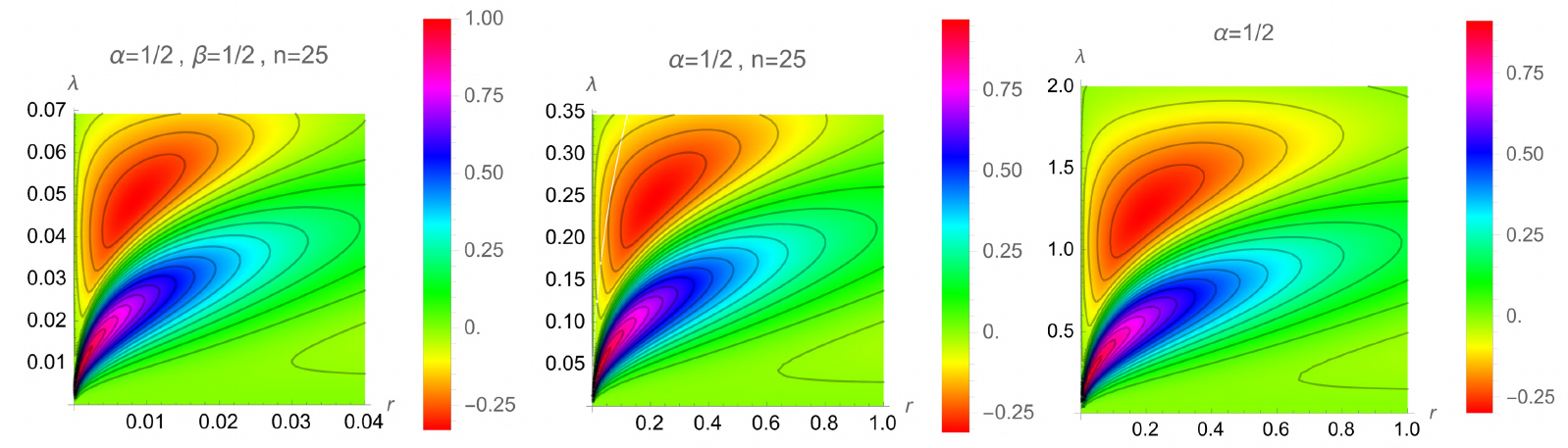}
        \captionsetup{justification=centering}
        \caption{\emph{(Left)} Plot of $(r,\lambda)\mapsto 2 \lambda \mathrm{cov}(r;\lambda^2)$, for the Jacobi ensemble with parameter $n=25$, $\alpha=1/2$, $\beta=1/2$.\\
        \emph{(Middle)} Plot of $(r,\lambda)\mapsto n^{3/2} 2 \lambda \mathrm{cov}(r;\lambda^2)$ for the Laguerre ensemble with parameter $n=25$, $\alpha=1/2$.\\ \emph{(Right)} Plot of $(r,\lambda)\mapsto 2 \lambda \mathrm{cov}^{(\infty)}(r;\lambda^2)$, for $w^{(\infty)}(x)=x^\alpha e^{-x}$, $\alpha=1/2$.}   
\end{figure}\label{fig: limiting cov}

\section*{Acknowledgement}
I thank Mario Kieburg and Arno Kuijlaars for their advice and feedback, as well as Sampad Lahiry and Mathieu Yahiaoui for fruitful discussions. This research is supported by the International Research Training Group (IRTG) between the University of Melbourne and KU Leuven and Melbourne Research Scholarship of University of Melbourne.

\appendix
\section{Consistency check of assumptions}\label{Appendix}

Let us show that, for P\'olya ensembles verifying \autoref{assum:exist HE} and for which the polynomials $\{p_j\}_{j=0,\ldots,n-1}$ are given by
\Eq{eq: polynomial polya}{
p_j\left(\frac{x }{\nu_n} \right)=\sum_{c=0}^{j} \binom{j}{c} \frac{(-x)^c}{\xi_n\nu_n^{c}\mathcal{M}w_n(c+1)},
}
the following lemma holds.
\lem{lem: asymp pj}{Let $n\in \mathbb{N}$, $n>2$, $x\in \mathbb{R}$, $R>0$. For a P\'olya ensemble verifying \autoref{assum:exist HE} with polynomials $\{p_j\}_{j=0,\ldots,n-1}$ given by \eqref{eq: polynomial polya}, we have 
\Eq{eq: integral pj limit}{
\lim_{R\to \infty} \lim_{n\to \infty}\int_{R}^{\infty} dt  \frac{ (1+t)}{\left(1+\frac{t}{n}\right)^{n+2}}\frac{1}{n}\sum_{j=0}^{n-1} \abs{p_j\left(\frac{x t}{n\nu_n} \right) }=0.
}
}
\begin{proof}[Proof of \autoref{lem: asymp pj}]
The idea of the following computation is to rewrite the sum of polynomials as another polynomial and use the \autoref{assum:exist HE} to show the non-vanishing contribution of the integrand in \eqref{eq: integral pj limit} remains in a bounded interval close to the origin as $n\to\infty$. The assumption \eqref{eq: lower bound h} excludes a fixed number of the first summands. This is not an issue as one can split the sum in two parts and deal with this fixed number of summands separately. Then, one can compute the $t$ integral explicitly in terms of an incomplete beta functions and treat the different regimes depending of the indices $j$. \\

Let us start by denoting
\Eq{}{
\alpha_{c,n}:=\xi_n\nu_n^{c}\mathcal{M}w_n(c+1).
}
We then use the summation identity
\Eq{}{
\sum_{j=0}^{n-1}\sum_{k=0}^j a_{k,j}=\sum_{j=0}^{n-1}\sum_{c=j}^{n-1} a_{j,c}
}
to express the sum of polynomials as another polynomial
\Eq{}{
\frac{1}{n}\sum_{j=0}^{n-1} \abs{p_j\left(\frac{x }{n\nu_n} \right) }\leq\frac{1}{n}\sum_{j=0}^{n-1}\sum_{k=0}^j \binom{j}{k} \frac{\abs{x}^k}{n^k\abs{\alpha_{k,n}}}=\sum_{j=0}^{n-1}\left(\frac{1}{n}\sum_{c=j}^{n-1} \binom{c}{j}\right) \frac{\abs{x}^j}{n^j \abs{\alpha_{j,n}}}.
}
Using the following bound,
\Eq{}{
\frac{1}{n}\sum_{c=j}^{n-1} \binom{c}{j}=\binom{n-1}{j}\frac{1}{j+1}\leq \binom{n-1}{j},
}
one gets to
\Eq{}{
\frac{1}{n}\sum_{j=0}^{n-1} \abs{p_j\left(\frac{x t}{n\nu_n} \right) }\leq \sum_{j=0}^{n-1} \binom{n-1}{j}\frac{\abs{xt}^j}{n^j \abs{\alpha_{j,n}}}.
}
Turning to the integral, we then get the following bound after proceeding with the change of variable $t\mapsto n t$,
\Eq{}{
\int_{R}^{\infty} dt  \frac{ (1+t)}{\left(1+\frac{t}{n}\right)^{n+2}}\frac{1}{n}\sum_{j=0}^{n-1} \abs{p_j\left(\frac{x t}{n\nu_n} \right) }\leq S_n(R)+T_n(R),
}
with
\Eq{eq: S_n}{
S_n(R):=2n^2\sum_{j=0}^{\lfloor R\rfloor-1} \int_{\frac{R}{n}}^{\infty} dt  \frac{ t }{\left(1+t\right)^{n+2}}\binom{n-1}{j}\frac{\abs{xt}^j}{ \abs{\alpha_{j,n}}},
}
and
\Eq{eq: Tn}{
T_n(R):=2n^2\sum_{j=\lfloor R\rfloor}^{n-1} \int_{\frac{R}{n}}^{\infty} dt  \frac{ t }{\left(1+t\right)^{n+2}}\binom{n-1}{j}\frac{\abs{xt}^j}{ h(j)}.
}
We have introduced, here, a cut-off at $j=\lfloor R\rfloor$, which serves to distinguish between two different asymptotic regimes. Note that we have used assumption \eqref{eq: lower bound h} to get \eqref{eq: Tn}. Indeed, as $j_*$ is fixed one can always choose $R$ big enough such that $R>j_*$.
\subsubsection*{Computation of the $t$ integral}
The $t$ integral can be computed term by term and yields an incomplete beta function
\Eq{eq: integ beta}{
\int_{\frac{R}{n}}^{\infty} dt\frac{ t^{j+1}}{\left(1+t\right)^{n+2}}= \mathrm{B}\left(\frac{n}{n+R};n-j,j+2\right)=\int_{0}^{\frac{n}{n+R}} dt\, t^{n-j-1}(1-t)^{j+1}.
}
To extract the $R$ dependence which will be crucial in the coming analysis, one needs to distinguish two different regimes as $j$ can grow with $n$. Hence the arbitrary choice of the cutoff $j=\lfloor R\rfloor$.

The first regime is for $j\in \llbracket 0, \lfloor R\rfloor-1\rrbracket$. As $j$ does not grow with $n$ here, one can use the asymptotic expansion of the beta function, given in \cite{Lopez1999},
\Eq{eq: bound beta fct}{
\mathrm{B}\left(\frac{n}{n+R};n-j,j+2\right)\underset{n\to\infty}{=}\frac{1}{n-j+1}\left(\frac{n}{n+R}\right)^{n-j} \left(\frac{R}{n+R}\right)^{j+1}\left[1+O\left(\frac{(j+1)!}{n-j+1}\right)\right] .
}
For the second regime, $j\in \llbracket \lfloor R\rfloor-1 ,n-1 \rrbracket$, one can simply use the identity
\Eq{}{
\binom{n-1}{j}=\frac{j+1}{n(n+1)}\frac{1}{\mathrm{B}(n-j,j+2)}
}
to get the following bound
\Eq{eq: bound beta fct 2}{
\mathrm{B}\left(\frac{n}{n+R};n-j,j+2\right)\leq \mathrm{B}(n-j,j+2)=\frac{j+1}{n(n+1)}\frac{1}{\binom{n-1}{j}}.
}

\subsubsection*{Computation of $S_n(R)$}
Let us first focus on $S_n(R)$ \eqref{eq: S_n}. The integral is given by the beta function \eqref{eq: integ beta} and one can split the sum in two parts to use assumption \eqref{eq: lower bound h} and get the bound
\Eq{}{
S_n(R)\leq S_n^{(1)}(R)+S_n^{(2)}(R),
}
with
\Eq{}{
S_n^{(1)}(R)&:=2n^2\sum_{j=0}^{j_*-1}\binom{n-1}{j}\frac{\abs{x}^j}{ \abs{\alpha_{j,n}}}\mathrm{B}\left(\frac{n}{n+R};n-j,j+2\right),\\
S_n^{(2)}(R)&:=2n^2\sum_{j=j_*}^{\lfloor R\rfloor-1}\binom{n-1}{j}\frac{\abs{x}^j}{ h(j)}\mathrm{B}\left(\frac{n}{n+R};n-j,j+2\right).
}
As all indices are in the interval $[0,R-1 ]$ one can use the expansion \eqref{eq: bound beta fct}. Combining it with the fact
\Eq{}{
\binom{n-1}{j}\frac{1}{n^j}\leq\frac{1}{j!},
}
along with the bound assumption \eqref{eq: bound Tildew} on $S_n^{(1)}(R)$, one gets, for all $n$ large enough,
\Eq{}{
S_n^{(1)}(R)\leq C_1 R\left(1-\frac{R}{n+R}\right)^{n-j_*+1} \sum_{j=0}^{j_*-1}\frac{(R\abs{x})^j}{ j!},
}
\Eq{}{
S_n^{(2)}(R)\leq C_2 R\left(1-\frac{R}{n+R}\right)^{n-\lfloor R\rfloor+1} \sum_{j=j_*}^{\lfloor R\rfloor-1}\frac{(R\abs{x})^j}{ j!\, h(j)},
}
where $C_1,C_2$ are some positive constants independent of $n$ and $R$. As $j_*$ is $n$-independent, taking the limit $n\to\infty$ on both sides of the two inequalities above yields
\Eq{eq: bound Sn1 and Sn2}{
\lim_{n\to\infty}S_n^{(1)}(R)+S_n^{(2)}(R)\leq C_1 R e^{-R}\sum_{j=0}^{j_*-1}\frac{(R\abs{x})^j}{ j!}+C_2 Re^{-R} \sum_{j=j_*}^{\lfloor R\rfloor-1}\frac{(R\abs{x})^j}{ j!\,h(j)}.
}
Let us show the RHS of the above inequality vanishes in the limit $R\to\infty$. Making use of Stirling approximation to get a rough lower bound of the factorial,
\Eq{}{
  \left(\frac{j}{e}\right)^j\leq \sqrt{2\pi j}\left(\frac{j}{e}\right)^j\leq \Gamma(j+1),\quad  j\geq 1,
}
one gets the following inequalities
\Eq{eq: inequalities bound}{
e^{-R} \sum_{j=j_*}^{\lfloor R\rfloor-1}\frac{(R\abs{x})^j}{ j!\,h(j)}\leq e^{-R/2} \sum_{j=j_*}^{\lfloor R\rfloor-1} \left(\frac{R}{ j}\exp(-\frac{R}{2j})\abs{x}e\right)^j\frac{1}{h(j)}\leq e^{-R/2} \sum_{j=j_*}^{\infty} \frac{\left(2\abs{x}\right)^j}{h(j)},
}
where we have used the fact
\Eq{}{
y\exp(-y/2)\leq 2e^{-1},\quad \forall y\geq 0.
}
Note that the RHS of \eqref{eq: inequalities bound} is indeed finite due to the root asymptotic of $h$ \eqref{eq: root asymp bound h}, which guarantees infinite radius of convergence of the series, by Cauchy-Hadamard theorem. Moreover, as $j_*$ is independent of $R$, the series itself is independent of $R$ and we have 
\Eq{}{
\lim_{R\to\infty}\left(C_1 R e^{-R}\sum_{j=0}^{j_*-1}\frac{(R\abs{x})^j}{ j!}+e^{-R/2} \sum_{j=j_*}^{\infty} \frac{\left(2\abs{x}\right)^j}{h(j)}\right)=0,
}
which finishes to show $\lim_{R\to\infty}\lim_{n\to\infty}S_n(R)=0$.
\subsubsection*{Computation of $T_n(R)$}
Turning to $T_n(R)$ \eqref{eq: Tn}, in this case, all the indices are in the interval $[R,n-1]$. Therefore, one can use the bound \eqref{eq: bound beta fct 2} to get
\Eq{}{
T_n(R)\leq 2\sum_{j=\lfloor R\rfloor}^{n-1} (j+1)\frac{\abs{x}^j}{ h(j)}.
}
Due to the root asymptotic of $h$ \eqref{eq: root asymp bound h}, by Cauchy-Hadamard theorem, the RHS is a convergent series as $n\to\infty$. Hence,
\Eq{}{
\lim_{n\to\infty} T_n(R)\leq 2\sum_{j=\lfloor R\rfloor}^{\infty} (j+1)\frac{\abs{x}^j}{ h(j)}.
}
The RHS of the above inequality being the remainder of a convergent series, taking the limit $R\to\infty$ finishes to show $\lim_{R\to\infty}\lim_{n\to\infty} T_n(R)=0$.

\end{proof}

\bibliography{main}

\end{document}